\documentclass{agtart_a}
\pdfoutput=1

\usepackage{pinlabel}
\usepackage[all]{xy}

\title[Surgery along Brunnian links]{On surgery along Brunnian links in $3$--manifolds}
\author[J-B Meilhan]{Jean-Baptiste Meilhan}
\givenname{Jean-Baptiste}
\surname{Meilhan}
\address{Research Institute for Mathematical Sciences\\Kyoto
University\\\newline
Kyoto 606-8502\\Japan}
\email{meilhan@kurims.kyoto-u.ac.jp}
\urladdr{}

\volumenumber{6}
\issuenumber{}
\publicationyear{2006}
\papernumber{84}
\startpage{2417}
\endpage{2453}

\doi{}
\MR{}
\Zbl{}

\keyword{3-manifolds}
\keyword{finite type invariants}
\keyword{Brunnian links}
\keyword{Goussarov-Vassiliev invariants}
\keyword{claspers}
\subject{primary}{msc2000}{57N10}
\subject{secondary}{msc2000}{57M27}

\received{30 May 2006}
\revised{}
\accepted{14 November 2006}
\published{13 December 2006}
\publishedonline{13 December 2006}
\proposed{}
\seconded{}
\corresponding{}
\editor{CPR}
\version{}

\arxivreference{math.GT/0603421}

\renewcommand{\subsubsection}[1]{\par\smallskip\stepcounter{subsubsection}\noindent{\bf\thesubsubsection\qua#1}\qua}

\let\xysavmatrix\xymatrix
\def\xymatrix{\disablesubscriptcorrection\xysavmatrix}
\AtBeginDocument{\let\bar\wbar\let\tilde\wtilde\let\hat\what}
\def\special#1{\relax}

\makeatletter
\def\cnewtheorem#1[#2]#3{\newtheorem{#1}{#3}[section]
\expandafter\let\csname c@#1\endcsname\c@thm}
\makeatother

\theoremstyle{plain}
\newtheorem{thm}{Theorem}[section]
\cnewtheorem{lem}[thm]{Lemma}
\cnewtheorem{prop}[thm]{Proposition}
\cnewtheorem{cor}[thm]{Corollary}
\theoremstyle{definition}
\cnewtheorem{defi}[thm]{Definition}
\theoremstyle{remark}
\cnewtheorem{rem}[thm]{Remark}
\cnewtheorem{claim}[thm]{Claim}
\cnewtheorem{summary}[thm]{Summary}
\cnewtheorem{example}[thm]{Example}
\newtheorem*{acknowledgments}{Acknowledgments}
\cnewtheorem{conv}[thm]{Convention}
\cnewtheorem{nota}[thm]{Notation}
\numberwithin{equation}{section}
\numberwithin{figure}{section}
\newcommand{\bd}{\begin{description}}
\newcommand{\ed}{\end{description}}
\newcommand{\ba}{\begin{array}}      \newcommand{\ea}{\end{array}}
\newcommand{\bc}{\begin{center}}     \newcommand{\ec}{\end{center}}
\newcommand{\be}{\begin{enumerate}}  \newcommand{\ee}{\end{enumerate}}
\newcommand{\beq}{\begin{eqnarray}}  \newcommand{\eeq}{\end{eqnarray}}
\newcommand{\beQ}{\begin{eqnarray*}} \newcommand{\eeQ}{\end{eqnarray*}}
\newcommand{\beqn}{\begin{equation}} \newcommand{\eeqn}{\end{equation}}
\newcommand{\bi}{\begin{itemize}}    \newcommand{\ei}{\end{itemize}}
\newcommand{\ov}{\overline}

\newcommand{\oz}{\otimes \mathbf{Z}[1 / 2]}

\newcommand{\oq}{\otimes \mathbf{Q}}

\makeop{Br}
\newcommand\Co{\operatorname{Co}}

\begin{document}

\begin{asciiabstract}
We consider surgery moves along (n+1)-component Brunnian links in
compact connected oriented 3-manifolds, where the framing of the each
component is 1/k for k in Z.  We show that no finite type invariant of
degree < 2n-2 can detect such a surgery move.  The case of two
link-homotopic Brunnian links is also considered.  We relate finite
type invariants of integral homology spheres obtained by such
operations to Goussarov-Vassiliev invariants of Brunnian links.
\end{asciiabstract}

\begin{htmlabstract}
We consider surgery moves along (n+1)&ndash;component Brunnian links in
compact connected oriented 3&ndash;manifolds, where the framing of the
components is in {1/k; k&isin;<b>Z</b>}.  We
show that no finite type invariant of degree < 2n-2 can detect such
a surgery move.  The case of two link-homotopic Brunnian links is also
considered.  We relate finite type invariants of integral homology
spheres obtained by such operations to Goussarov&ndash;Vassiliev invariants
of Brunnian links.
\end{htmlabstract}

\begin{abstract}
We consider surgery moves along $(n+1)$--component Brunnian links in
compact connected oriented $3$--manifolds, where the framing of the
components is in $\{ \frac{1}{k}\textrm{ ; } k\in \mathbf{Z} \}$.  We
show that no finite type invariant of degree $< 2n-2$ can detect such
a surgery move.  The case of two link-homotopic Brunnian links is also
considered.  We relate finite type invariants of integral homology
spheres obtained by such operations to Goussarov--Vassiliev invariants
of Brunnian links.
\end{abstract}
\maketitle

\section{Introduction}\label{sec1}
In \cite{O}, Ohtsuki introduced the notion of finite type invariants of integral homology spheres as an attempt to unify 
the topological invariants of these objects, in the same way as Goussarov--Vassiliev invariants provide a unified point of view 
on invariants of knots and links.  
This theory was later generalized to all oriented $3$--manifolds by Cochran and Melvin \cite{CM}.  

Goussarov and Habiro developed independently another finite type invariants theory for compact connected oriented $3$--manifolds, 
which essentially coincides with the Ohtsuki theory in the case of integral homology spheres  
\cite{GGP,goussarov,habiro}.  This theory comes equipped with a new and powerful tool 
called \emph{calculus of clasper}, which uses embedded graphs carrying some surgery instruction.   
Surgery moves along claspers 
define a family of (finer and finer) equivalence relations among $3$--manifolds, called 
\emph{$Y_k$--equivalence}, which gives a good idea of the information contained by finite type invariants: 
two compact connected oriented $3$--manifolds are not distinguished by invariants of degree $<k$ if they are 
$Y_{k}$--equivalent \cite{G2,habiro}.  These two conditions become equivalent when dealing with integral homology spheres.  

Recall that a link $L$ is \emph{Brunnian} if any proper sublink of $L$ is trivial.  
In some sense, an $n$--component Brunnian link is a `pure $n$--component linking'.
In this paper we consider those compact connected oriented $3$--manifolds which are obtained by 
surgery along a Brunnian link.  
For a fixed number of components, we study which finite type invariants 
(ie of which degree) can vary under such an operation.  

Let $m=(m_1,...,m_n)\in \mathbf{Z}^{n}$ be a collection of $n$ integers.  
Given a null-homologous, ordered $n$--component link $L$ in a compact connected oriented $3$--manifold $M$, 
denote by $(L,m)$ the link $L$ with framing $\frac{1}{m_i}$ on the $i^{th}$ component ; $1\le i\le n$.  
We denote by $M_{(L,m)}$ the $3$--manifold obtained from $M$ by surgery along the framed link $(L,m)$.  
We say that $M_{(L,m)}$ is obtained from $M$ by \emph{$\frac{1}{m}$--surgery} along the link $L$.   
\begin{thm} \label{2n-2}
  Let $n\ge 2$ and $m\in \mathbf{Z}^{n+1}$.  
  Let $L$ be an $(n+1)$--component Brunnian link in a compact, connected, oriented $3$--manifold $M$.  

  For $n=2$, $M_{(L,m)}$ and $M$ are $Y_{1}$--equivalent.  

  For $n\ge 3$, $M_{(L,m)}$ and $M$ are $Y_{2n-2}$--equivalent.  
  Consequently, they cannot be distinguished by any finite type invariant of degree $< 2n-2$.  
\end{thm}
\noindent Note that, for any Brunnian link $L$ in $M$, we have $M_{(L,m)}\cong M$ if $m_i=0$ for some $1\le i\le n+1$.  
In this case, the statement is thus vacuous.  

Two links are \emph{link-homotopic} if they are related by a sequence of isotopies and 
self-crossing changes, ie, crossing changes involving two strands of the same component.  
We obtain the following.  
\begin{thm} \label{2n-1}
  Let $n\ge 2$ and $m\in \mathbf{Z}^{n+1}$.  
  Let $L$ and $L'$ be two link-homotopic $(n+1)$--component Brunnian links in a compact, connected, oriented $3$--manifold $M$.    
  Then $M_{(L,m)}$ and $M_{(L',m)}$ are $Y_{2n-1}$--equivalent.  
  Consequently, they cannot be distinguished by any finite type invariant of degree $< 2n-1$.  
\end{thm}
\noindent Actually, for integral homology spheres, the theorem is still true when ``$2n-1$'' is replaced by ``$2n$''.  
(It follows from the last observation of \fullref{mariethomas}.)  

In the latter part of the paper, we study the relation 
between the above results and Goussarov--Vassiliev invariants of Brunnian links.  

Let $\mathbf{Z}\mathcal{L}(n)$ be the free $\mathbf{Z}$--module generated by the set of isotopy classes of $n$--component links 
in $S^3$.  
The theory of Goussarov--Vassiliev invariants of links involves a descending filtration 
  $$ \mathbf{Z}\mathcal{L}(n)=J_0(n)\supset J_1(n)\supset J_2(n)\supset ... $$
called \emph{Goussarov--Vassiliev filtration} (see \fullref{brunnianpart}).   
In a previous paper, Habiro and the author introduced the so-called \emph{Brunnian part} $\Br(\overline{J}_{2n}(n+1))$ 
of $J_{2n}(n+1)/J_{2n+1}(n+1)$, which is defined as the $\mathbf{Z}$--submodule generated by elements $[L-U]_{J_{2n+1}}$ 
where $L$ is an $(n+1)$--component Brunnian link and $U$ is the $(n+1)$--component unlink \cite{hm2}.  
Further, we constructed a linear map
  $$ h_n\co  \mathcal{A}^c_{n-1}(\emptyset)\longrightarrow \Br(\overline{J}_{2n}(n+1)), $$
where $\mathcal{A}^c_{n-1}(\emptyset)$ is a $\mathbf{Z}$--module of connected trivalent
diagrams with $2n-2$ vertices.  $h_n$ is an isomorphism over $\mathbf{Q}$ for $n\ge 2$.  
See \fullref{GV} for precise definitions.  

Let $\overline{\mathcal{S}}_{k}$ be the abelian group of $Y_{k+1}$--equivalence classes of integral homology spheres which are 
$Y_{k}$--equivalent to $S^3$.  $\overline{\mathcal{S}}_{2k+1}=0$ for all $k\ge 1$, and it is well known that 
$\overline{\mathcal{S}}_{2k}$ is isomorphic to $\mathcal{A}^c_{k}(\emptyset)$ when tensoring by $\mathbf{Q}$. 
See \fullref{sec:Sk}.  
There is therefore an isomorphism over $\mathbf{Q}$ from $\Br(\overline{J}_{2n}(n+1))$ to $\overline{\mathcal{S}}_{2n-2}$, 
for $n\ge 2$.  The next theorem states that this isomorphism is induced by $(+1)$--framed surgery. 

For a null-homologous ordered link $L$ in a compact connected oriented $3$--manifold $M$, 
denote by $(L,+1)$ the link $L$ with all components having framing $+1$.  
\begin{thm} \label{homspheres}
  For $n\ge 2$, the assignment 
   $$ [L-U]_{J_{2n+1}} \mapsto [S^3_{(L,+1)}]_{Y_{2n-1}} $$
  defines an isomorphism 
    $$ \kappa_n\co   \Br(\overline{J}_{2n}(n+1))\oq\longrightarrow \overline{\mathcal{S}}_{2n-2}\oq. $$
\end{thm}
We actually show that these two $\mathbf{Q}$--modules are isomorphic to the so-called `connected part' of the Ohtsuki filtration, 
by using the abelian group $\mathcal{A}^c_{n-1}(\emptyset)$.  See \fullref{FTIZHS} for definitions and statements.  

The rest of this paper is organized as follows.  

In \fullref{claspers}, we give a brief review of the theory of claspers, both for compact connected oriented $3$--manifolds 
and for links in a fixed manifold.  
In \fullref{secspecial}, we study the $Y_k$--equivalence class of integral homology spheres obtained by surgery along claspers with several 
\emph{special leaves}.  This section can be read separately from the rest of the paper and might be of independent interest.  
In \fullref{+1surgery}, we use the main result of section 3 to prove Theorems \ref{2n-2} and \ref{2n-1}.  
In \fullref{GV}, we recall several results obtained by Habiro and the author in \cite{hm2}.  
In \fullref{FTIZHS}, we define the material announced above and prove \fullref{homspheres}.  
In \fullref{zeproof}, we give the (technical) proof of \fullref{linear2n+1}.  
\begin{acknowledgments}
  The author is grateful to 
  Kazuo Habiro 
  for many helpful conversations and comments on an early version of this paper.  
  He was supported by a Postdoctoral
  Fellowship and a Grant-in-Aid for Scientific Research of the Japan
  Society for the Promotion of Science.
\end{acknowledgments}
\section{Claspers} \label{claspers}
Throughout this paper, all $3$--manifolds will be supposed to be compact, connected and oriented. 

\subsection{Clasper theory for $3$--manifolds} 
Let us briefly recall from \cite{GGP,goussarov,habiro} the fundamental notions of clasper theory for $3$--manifolds.  
\begin{defi}
A {\em clasper} in a $3$--manifold $M$ is an embedding 
$$ G\co  F\longrightarrow \operatorname{int} M $$
of a compact (possibly unorientable) surface $F$.  $F$ is decomposed into \emph{constituents} connected by disjoint bands called 
\emph{edges}.  
Constituents are disjoint connected subsurfaces, either annuli or disks: 
\begin{itemize}
\item A {\em leaf} is an annulus with one edge attached.
\item A {\em node} is a disk with three edges attached.
\item A {\em box} is a disk with \emph{at least} three edges attached, one being 
  distinguished with the others.  This distinction is done by drawing a box as a rectangle.  
\end{itemize}
\end{defi}
\noindent Observe that this definition slightly extends the one in \cite{habiro}, where a box has always three edges attached.  

We will make use of the drawing convention for claspers of \cite[Figure 7]{habiro}, except for the following: 
a $\oplus$ (resp.\ $\ominus$) on an edge represents a positive (resp.\ negative) half-twist.  This replaces the 
convention of a circled $S$ (resp.\ $S^{-1}$) used in \cite{habiro}.  
\subsubsection{Surgery along claspers}
\begin{figure}[ht!]
  \bc
\labellist\small
\pinlabel (a) [t] at 46 26
\pinlabel (b) [t] at 180 26
\pinlabel (c) [t] at 310 26
\endlabellist
    \includegraphics{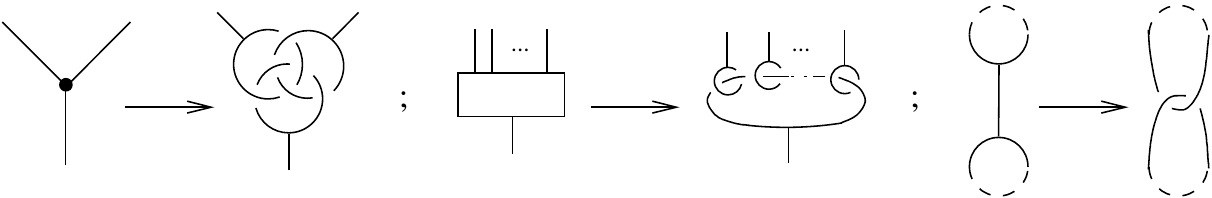}
    \caption{Constructing the framed link associated to a clasper} \label{link}
  \ec
  \end{figure}

Given a clasper $G$ in $M$, we can construct, in a regular 
neighborhood of the clasper, an associated framed link $L_G$ as follows.  
First, replace each node and box of $G$ by leaves as shown in \fullref{link} (a) and (b).  We obtain 
a union of $\textsf{I}$--shaped claspers, one for each edge of $G$.  $L_G$ is obtained by replacing each of these 
$\textsf{I}$--shaped claspers by a $2$--component framed link as shown in \fullref{link} (c).\footnote{Here 
  and throughout the paper, blackboard framing convention is used.} 
  
\emph{Surgery along the clasper $G$} is defined to be surgery along $L_G$.  

In \cite[Proposition 2.7]{habiro}, Habiro gives a list of 12 moves on claspers which gives \emph{equivalent} claspers, that is 
claspers with diffeomorphic surgery effect.  
We will freely use \emph{Habiro's moves} (which are essentially derived from Kirby calculus) 
by referring to their numbering in Habiro's paper.   
\subsubsection{The $Y_k$--equivalence}
For $n\ge 1$, a {\em $Y_n$--graph} is a connected clasper $G$ without boxes and with $n$ nodes, where 
a connected clasper is a clasper whose underlying surface is connected.  
The integer $n$ is called the \emph{degree} of $G$.  

\begin{figure}[ht!]
  \bc
\labellist\small
\pinlabel $f$ [t] at 0 14
\pinlabel $f'$ [t] at 112 14
\pinlabel $T$ [b] at 68 20
\endlabellist
    \includegraphics[scale=1.2]{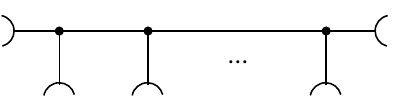}
    \caption{A linear tree $T$ and its two ends $f$ and $f'$} \label{LINE}
  \ec
  \end{figure}

A \emph{$Y_k$--tree} is a $Y_k$--graph $T$ such that the union of edges and nodes of $T$ is simply connected.  
For $k\ge 3$, we say that a $Y_k$--tree $T$ in a $3$--manifold $M$ is \emph{linear} if there is a $3$--ball in $M$ 
which intersects the edges and nodes of $T$ as shown in \fullref{LINE}.  The leaves denoted by $f$ and $f'$ in the figure 
are called the \emph{ends} of $T$.

A \emph{$Y_k$--forest} is a clasper $T=T_1\cup ...\cup T_p$ ($p\ge 0$), where $T_i$ is a $Y_{k_i}$--tree ($p\ge i\ge 1$), 
such that $min_{1\le i\le p} k_i=k$. 

A \emph{$Y_k$--subtree} $T$ of a clasper $G$ is a connected union of leaves, nodes and edges of $G$ 
such that the union of edges and nodes of $T$ is simply connected and such that $T$ intersects $\overline{G\setminus T}$ along the 
attaching region of some edges of $T$, called \emph{branches}.  

A surgery move on $M$ along a $Y_k$--graph $G$ is called a \emph{$Y_k$--move}.  
For example, a $Y_1$--move is equivalent to Matveev's Borromean surgery \cite{matv}.   

The \emph{$Y_k$--equivalence} is the equivalence relation on $3$--manifolds 
generated by $Y_k$--moves and orientation-preserving diffeomorphisms.  
This equivalence relation becomes finer as $k$ increases: if $k\le l$ and if $M\sim_{Y_l}N$, then we also have $M\sim_{Y_k}N$.  

Recall that `trees do suffice to define the $Y_k$--equivalence'.  That is, $M\sim_{Y_k}N$ implies that there exists 
a $Y_k$--forest $F$ in $M$ such that $M_F\cong N$.  
\subsection{Clasper theory for links} 
Another aspect of the theory of claspers is that it allows to study links in a \emph{fixed} manifold.  
For this we use a slightly different type of claspers.  
\begin{defi} \label{fini}
Let $L$ be a link in a $3$--manifold $M$, and let $G$ be a clasper in $M$ which is disjoint from $L$. 
A \emph{disk-leaf} of $G$ is 
a leaf $l$ of $G$ which is an unknot bounding a disk $D$ in $M$ with respect to which it is $0$--framed.\footnote{ 
Here we regard a leaf, which is an embedded annulus, as a knot with a framing.  }
We call $D$ the \emph{bounding disk} of $f$.  
The interior of $D$ is disjoint from $G$ and from any other bounding disk, 
but it may intersect $L$ transversely. 
For convenience, we say that a disk-leaf $f$ \emph{intersects} the link $L$ when its bounding disk does.  
\end{defi}

A \emph{$C_n$--tree} (resp.\ \emph{linear $C_n$--tree}) for a link $L$ in a 
$3$--manifold $M$ is a $Y_{n-1}$--tree (resp.\ \emph{linear $Y_{n-1}$--tree}) in $M$ 
such that each of its leaves is a disk-leaf.  

Given a $C_n$--tree $C$ in $M$, there exists a canonical diffeomorphism between $M$ and the manifold $M_C$.  
So surgery along a $C_n$--tree can be regarded as a local move on links in the manifold $M$.  

A $C_n$--tree $G$ for a link $L$ is {\em simple} (with respect to $L$) if each disk-leaf of $G$ intersects $L$ exactly once.  

A surgery move on a link $L$ along a $C_k$--tree is called a \emph{$C_k$--move}.  
The \emph{$C_k$--equivalence} is the equivalence relation on links generated by the $C_k$--moves and isotopies.  
As in the case of manifolds, the $C_n$--equivalence relation implies the $C_k$--equivalence if $1\le k\le n$.
For more details, see \cite{G2,habiro}.  
\subsection{Some technical lemmas}
In this subsection, we state several technical lemmas about claspers. 

First, we introduce several moves on claspers which produce equivalent claspers, like the 12 Habiro's moves.  
In each of the next three statements, the figure represents two claspers in a given $3$--manifold 
which are identical outside a $3$--ball, where they are as depicted.  
\begin{lem} \label{slide}
  The move of \fullref{move1}  produces equivalent claspers.  
  \begin{figure}[ht!]
  \bc
    \includegraphics{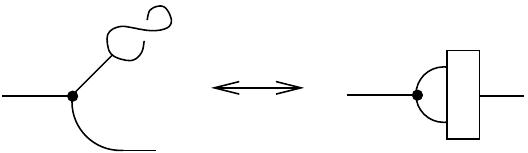}
    \caption{} \label{move1}
  \ec
  \end{figure}
\end{lem}
\noindent This is an immediate consequence of \cite[Theorem 3.1]{GGP} (taking into account that the convention used in \cite{GGP} for 
the definition of the surgery link associated to a clasper is the opposite of the one used in the present paper). 
\begin{lem} \label{move}
  The move of \fullref{move2} produces equivalent claspers.  
  \begin{figure}[ht!]
  \bc
   \includegraphics{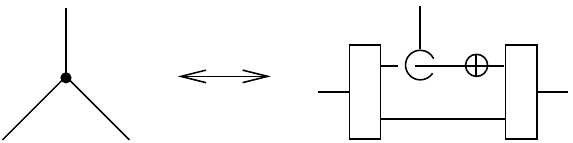} 
   \caption{} \label{move2}
  \ec
  \end{figure}
\end{lem}
\noindent This move is, in some sense, the inverse of Habiro's move 12.  
See also Figure 25 of \cite{CT}, where a similar move appears.  
\begin{proof}
 Consider the clasper on the right-hand side of \fullref{move2}.  By replacing the two boxes by leaves as shown in 
 \fullref{link} (b) and applying Habiro's move 1, we obtain the clasper depicted on the left-hand side of \fullref{proofmove}.  
   \begin{figure}[ht!]
   \bc
   \includegraphics{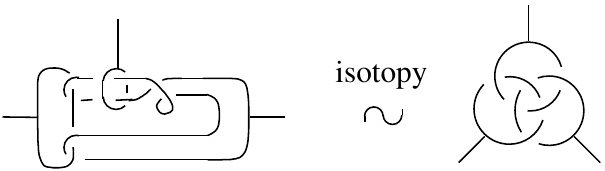} 
   \caption{} \label{proofmove}
   \ec
   \end{figure} 
 Now, the three leaves depicted in this figure form a $3$--component link which is isotopic to the Borromean link.  
 As shown in \fullref{link} (a), this is equivalent to a node.  
\end{proof}
\begin{lem} \label{lass}
  The moves of \fullref{assoc} produce equivalent claspers.  
   \begin{figure}[ht!]
   \bc
   \includegraphics{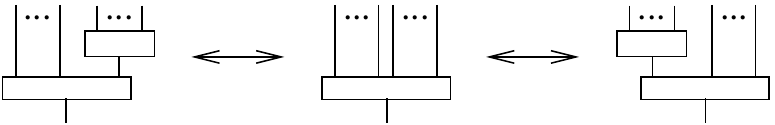} 
   \caption{The associativity of boxes} \label{assoc}
   \ec
   \end{figure} 
\end{lem} 
\noindent This `associativity' property of boxes is easily checked using \fullref{link} (b) -- see Figure 37 of \cite{habiro}.  

The next lemma deals with crossing change operations on claspers.  
A crossing change is a local move as illustrated in \fullref{fcc}.  
The proof is omitted, as it uses the same techniques as in \cite[Section 4]{habiro} (where similar statements appear).  
See also \cite[Section 1.4]{these}.  
 \begin{lem} \label{crossingchange}
  Let $T_1\cup T_2$ be a disjoint union of a $Y_{k_1}$--tree and a $Y_{k_2}$--tree in a $3$--manifold $M$.     
   \begin{figure}[ht!]
   \bc
\begin{picture}(0,0)%
\includegraphics{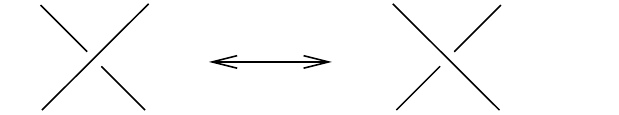}%
\end{picture}%
\setlength{\unitlength}{3947sp}%
\def\SetFigFont#1#2#3#4#5{\small}%
\begin{picture}(2419,534)(574,232)
\put(2210,253){\makebox(0,0)[lb]{\smash{\SetFigFont{12}{14.4}{\rmdefault}{\mddefault}{\updefault}
$T'_1$
}}}
\put(2993,254){\makebox(0,0)[lb]{\smash{\SetFigFont{12}{14.4}{\rmdefault}{\mddefault}{\updefault}
$T'_2$
}}}
\put(1287,258){\makebox(0,0)[lb]{\smash{\SetFigFont{12}{14.4}{\rmdefault}{\mddefault}{\updefault}
$T_2$
}}}
\put(574,258){\makebox(0,0)[lb]{\smash{\SetFigFont{12}{14.4}{\rmdefault}{\mddefault}{\updefault}
$T_1$
}}}
\end{picture}
\caption{A crossing change} \label{fcc}
   \ec
   \end{figure}
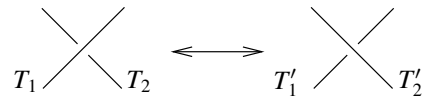 
  Let $T'_1\cup T'_2$ be obtained by a crossing change $c$ of an edge \emph{or} a leaf of 
  $T_1$ with an edge \emph{or} a leaf of $T_2$ (see \fullref{fcc}), 
  and let $C\in \{0,1,2\}$ denotes the number of edges involved in the crossing change $c$.  Then 
   \begin{enumerate}
     \item $M_{T_1\cup T_2}\sim_{Y_{k_1+k_2+C}} M_{T'_1\cup T'_2}$.
     \item $M_{T_1\cup T_2}\sim_{Y_{k_1+k_2+C+1}} M_{T'_1\cup T'_2\cup T}$, 
           where $T$ is a parallel copy, disjoint from $T'_1\cup T'_2$, of some $Y_{k_1+k_2+C}$--tree $\tilde{T}$ 
           obtained as follows: 
  	\begin{enumerate}
	     \item If $c$ involves an edge $e_1$ of $T_1$ and an edge $e_2$ of $T_2$, then $C=2$ and 
	           $\tilde{T}$ is obtained by inserting a node $n_1$ in $e_1$ and a node $n_2$ in $e_2$, 
	           and connecting $n_1$ and $n_2$ by an edge.
	     \item If $c$ involves an edge $e$ of $T_1$ and a leaf $f$ of $T_2$, then $C=1$ and 
	           $\tilde{T}$ is obtained by inserting a node $n$ in $e$, and connecting $n_1$ to the edge incident to $f$.  
	     \item If $c$ involves a leaf $f_1$ of $T_1$ and a leaf $f_2$ of $T_2$, then $C=0$ and 
	           $\tilde{T}$ is obtained by connecting the edges incident to $f_1$ and $f_2$.  
	\end{enumerate}
    \end{enumerate}
 \end{lem}
\begin{rem}
This lemma is only valid for trees.  
However, if we are given graphs or subtrees instead, observe that it suffices to use Habiro's move 2 to obtain 
equivalent trees.  So in this paper, whenever we apply \fullref{crossingchange} to graphs or subtrees, 
it implicitly means that we apply the lemma to some equivalent trees obtained by Habiro's move 2.  
\end{rem}
The next result follows from \fullref{crossingchange} and \cite[Proposition 2.7]{habiro}. See also \cite{GGP,ohtsuki}.  
\begin{lem} \label{lemtwist}
Let $G$ be a $Y_k$--tree in a $3$--manifold $M$, and let $G_{+}$ be a $Y_k$--tree obtained from $G$ 
by inserting a positive half twist in an edge.  Then 
  $$ M_{G\cup \tilde{G}_+}\sim_{Y_{k+1}} M, $$
where $\tilde{G}_+$ is obtained from $G_+$ by an isotopy so that it is disjoint from $G$.  
\end{lem}
\subsection{The IHX relation for $Y_k$--graphs} \label{subihx} 
We have the following version of the IHX relation for $Y_k$--graphs.  
\begin{lem} \label{ihx}
 Let $I$, $H$ and $X$ be three $Y_k$--graphs in a $3$--manifold $M$, 
 which are identical except in a $3$--ball where they look as depicted in \fullref{fig:ihx}.  
 Then 
 $$ M_I\sim_{Y_{k+1}} M_{H\cup \tilde{X}}, $$
 \noindent where $\tilde{X}$ is obtained from $X$ by an isotopy so that it is disjoint from $H$.  
   \begin{figure}[ht!]
   \bc
   \includegraphics{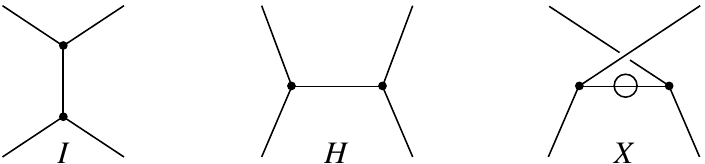} 
   \caption{The three $Y_k$--graphs $I$, $H$ and $X$} \label{fig:ihx}
   \ec
   \end{figure}
\end{lem}
\noindent Various similar statements appear in the literature.  
For example, an IHX relation is proved in \cite{GGP} at the level of finite type invariants, in \cite{CT} for 
$C_n$--trees (see also \cite{G2}), and in \cite[pages 397--398]{ohtsuki} for $Y_n$--graphs without leaves.  
\begin{proof}
 For simplicity, we give the proof for the case of $Y_2$--trees.  
 In the general case, the proof uses the same arguments as below, together with the zip construction 
 (\cite[Section 3]{habiro}, see also \cite[Section 4.2]{CT}).  

 Consider the $Y_2$--tree $I$, and apply \fullref{move} at one of its nodes.  Then, apply Habiro's move 11 
 so that we obtain the clasper $G_1\sim I$ depicted in \fullref{pihx1}.  
 By an isotopy and Habiro's move 7,  $G_1$ is seen to be equivalent to the clasper $G_2$ of \fullref{pihx1}.  
   \begin{figure}[ht!]
   \bc
\begin{picture}(0,0)%
\includegraphics{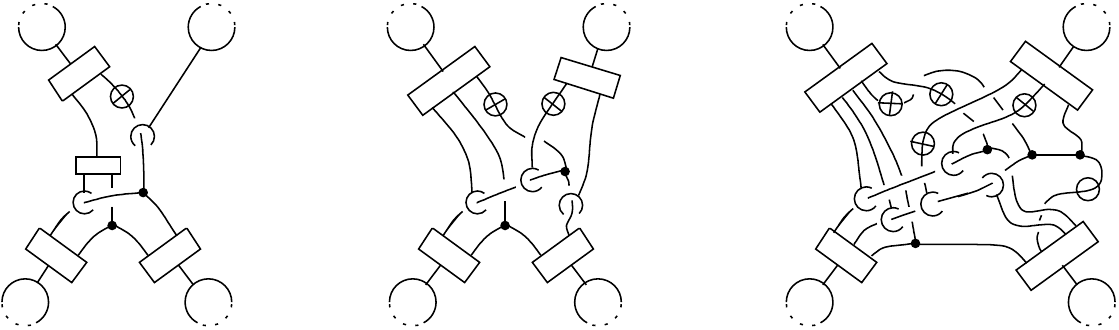}%
\end{picture}%
\setlength{\unitlength}{3947sp}%
\begingroup\makeatletter\ifx\SetFigFont\undefined%
\gdef\SetFigFont#1#2#3#4#5{%
  \reset@font\fontsize{#1}{#2pt}%
  \fontfamily{#3}\fontseries{#4}\fontshape{#5}%
  \selectfont}%
\fi\endgroup%
\begin{picture}(5362,1576)(502,-1796)
\put(1827,-1116){\makebox(0,0)[lb]{\smash{\SetFigFont{10}{12.0}{\rmdefault}{\mddefault}{\updefault}{\color[rgb]{0,0,0}$\sim$}%
}}}
\put(3732,-1116){\makebox(0,0)[lb]{\smash{\SetFigFont{10}{12.0}{\rmdefault}{\mddefault}{\updefault}{\color[rgb]{0,0,0}$\sim$}%
}}}
\put(3332,-1266){\makebox(0,0)[lb]{\smash{\SetFigFont{10}{12.0}{\rmdefault}{\mddefault}{\updefault}{\color[rgb]{0,0,0}$f$}%
}}}
\put(2832,-1751){\makebox(0,0)[lb]{\smash{\SetFigFont{10}{12.0}{\rmdefault}{\mddefault}{\updefault}{\color[rgb]{0,0,0}$G_2$}%
}}}
\put(952,-1746){\makebox(0,0)[lb]{\smash{\SetFigFont{10}{12.0}{\rmdefault}{\mddefault}{\updefault}{\color[rgb]{0,0,0}$G_1$}%
}}}
\put(5031,-1751){\makebox(0,0)[lb]{\smash{\SetFigFont{10}{12.0}{\rmdefault}{\mddefault}{\updefault}{\color[rgb]{0,0,0}$G_3$}%
}}}
\put(5737,-936){\makebox(0,0)[lb]{\smash{\SetFigFont{10}{12.0}{\rmdefault}{\mddefault}{\updefault}{\color[rgb]{0,0,0}$T_H$}%
}}}
\end{picture}
   \caption{} \label{pihx1}
   \ec
   \end{figure}
 Consider the leaf of $G_2$ denoted by $f$ in the figure.  
 By an application of Habiro's move 12 at $f$, followed by moves 7 and 11, 
 we obtain the clasper $G_3$ of \fullref{pihx1}.  Observe that $G_3$ contains a $Y_2$--subtree $T_{H}$.  
 By Habiro's move 6, \fullref{lemtwist} and \fullref{crossingchange} (1), we have 
  $$M_{G_3}\sim_{Y_3} M_{H\cup G_4},$$ 
 where $G_4$ is the clasper depicted in \fullref{pihx2}.  
   \begin{figure}[ht!]
   \bc
\begin{picture}(0,0)%
\includegraphics{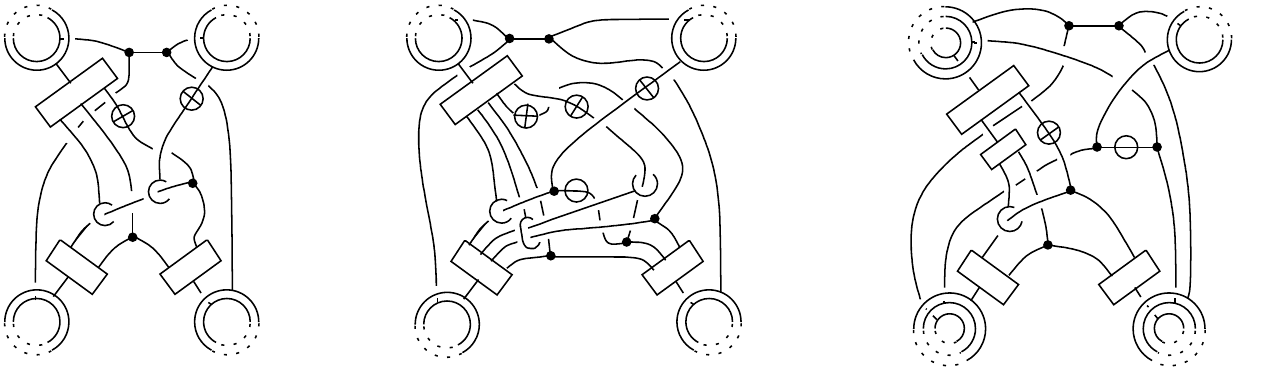}%
\end{picture}%
\setlength{\unitlength}{3947sp}%
\begingroup\makeatletter\ifx\SetFigFont\undefined%
\gdef\SetFigFont#1#2#3#4#5{%
  \reset@font\fontsize{#1}{#2pt}%
  \fontfamily{#3}\fontseries{#4}\fontshape{#5}%
  \selectfont}%
\fi\endgroup%
\begin{picture}(5906,1748)(-35,-958)
\put(708,-314){\makebox(0,0)[lb]{\smash{\SetFigFont{10}{12.0}{\rmdefault}{\mddefault}{\updefault}
$f'$
}}}
\put(1393,-154){\makebox(0,0)[lb]{\smash{\SetFigFont{10}{12.0}{\rmdefault}{\mddefault}{\updefault}
$\sim$
}}}
\put(3914,-98){\makebox(0,0)[lb]{\smash{\SetFigFont{10}{12.0}{\rmdefault}{\mddefault}{\updefault}
;
}}}
\put(4928,-584){\makebox(0,0)[lb]{\smash{\SetFigFont{10}{12.0}{\rmdefault}{\mddefault}{\updefault}
$G_0$
}}}
\put(5269,-136){\makebox(0,0)[lb]{\smash{\SetFigFont{10}{12.0}{\rmdefault}{\mddefault}{\updefault}
$\tilde{X}$
}}}
\put(1088,134){\makebox(0,0)[lb]{\smash{\SetFigFont{10}{12.0}{\rmdefault}{\mddefault}{\updefault}
$H$
}}}
\put(3368,134){\makebox(0,0)[lb]{\smash{\SetFigFont{10}{12.0}{\rmdefault}{\mddefault}{\updefault}
$H$
}}}
\put(5641,186){\makebox(0,0)[lb]{\smash{\SetFigFont{10}{12.0}{\rmdefault}{\mddefault}{\updefault}
$H$
}}}
\put(2693,-32){\makebox(0,0)[lb]{\smash{\SetFigFont{10}{12.0}{\rmdefault}{\mddefault}{\updefault}
$T_X$
}}}
\put(2620,-672){\makebox(0,0)[lb]{\smash{\SetFigFont{10}{12.0}{\rmdefault}{\mddefault}{\updefault}
$G_5$
}}}
\put(493,-699){\makebox(0,0)[lb]{\smash{\SetFigFont{10}{12.0}{\rmdefault}{\mddefault}{\updefault}
$G_4$
}}}
\end{picture}
   \caption{} \label{pihx2}
   \ec
   \end{figure}
 Now, consider the leaf $f'$ of $G_4$ (see the figure). Apply Habiro's move 12 at $f'$ and moves 7 and 11, just as we 
 did previously for the clasper $G_2$.  
 The resulting clasper $G_5\sim G_4$ contains a $Y_2$--subtree $T_{X}$ (see \fullref{pihx2}).  
 As above, we obtain by Lemmas \ref{lemtwist} and \ref{crossingchange} (1): 
   $$ M_{H\cup G_5}\sim_{Y_3} M_{H\cup \tilde{X}\cup G_0},$$  
 where $G_0$ is represented in the right-hand side of \fullref{pihx2}.  
 By Habiro's moves 11 and 4, we obtain that $M_{H\cup \tilde{X}\cup G_0}\cong M_{H\cup \tilde{X}}$.  
\end{proof}
One can check the following slightly stronger version of \fullref{ihx} when $I$, $H$ and $X$ are 
three $Y_k$--trees (Note that Habiro's move 2 always allows us to have this condition satisfied).  
\begin{lem} \label{corihx}
 Let $I$, $H$, $X$ and $\tilde{X}$ be four $Y_k$--trees in a $3$--manifold $M$ as in \fullref{ihx}. Then 
  $$ M_I\sim_{Y_{k+2}} M_{H\cup \tilde{X}\cup F}, $$
 where $F$ 
 is a union of 
 disjoint $Y_{k+1}$--trees.  
 Each $Y_{k+1}$--tree $T$ in $F$ is obtained from either $H$ or $X$ 
 by taking a parallel copy $f$ of one of its leaves, inserting a node $n$ in one of its edges,  
 connecting $n$ and $f$ by an edge, and performing an isotopy so that $T$ is disjoint from $H$, $\tilde{X}$ and $F\setminus T$.  
\end{lem}
\noindent
Consider for example the case of $Y_2$--trees, as in the proof of \fullref{ihx}.  
We saw there that $M_I\cong M_{G_3}\sim_{Y_3} M_{H\cup G_4}$, 
where $G_3$ and $G_4$ are depicted in \fullref{pihx1} and \ref{pihx2}.  
Observe that $H\cup G_4$ is obtained from $G_3$ by several Habiro's moves and three crossing changes between an edge 
of the $Y_2$--subtree $T_H$ and some leaf of $G_3$.  
So by (2) of \fullref{crossingchange} (and Habiro's move 5) one can check that 
 $$ M_{G_3}\sim_{Y_4} M_{H\cup G_4\cup F'}, $$
where $F'$ consists of three $Y_3$--trees obtained as described in the statement of the Lemma. 
For similar reasons, (2) of \fullref{crossingchange} implies that the clasper $G_5\sim G_4$ depicted in \fullref{pihx2} satisfies 
$ M_{H\cup G_5}\sim_{Y_4} M_{H\cup \tilde{X}\cup G_0\cup F''}, $ 
where $F''$ is a union of $Y_3$--trees of the desired form.  
This implies \fullref{corihx} for $k=2$.  
\section{Surgery along $Y_n$--trees with special leaves} \label{secspecial}
In this section, we study $3$--manifolds obtained by surgery along $Y_n$--trees containing a particular type of leaves.  
\subsection{$m$--special leaves} \label{leaves}
Suppose we are given a clasper $G$ in a $3$--manifold $M$.  
\begin{defi}
Let $m\in \mathbf{Z}$.  
An \emph{$m$--special leaf with respect to $G$} is a leaf $f$ of $G$ which is an 
unknot bounding a disk $D$ in $M$ with respect to which it is $m$--framed,\footnote{ 
Here, as in \fullref{fini}, we regard a leaf as a knot with a framing.  }
such that the interior of $D$ is disjoint from $G\setminus f$.    
$D$ is called the \emph{bounding disk} of $f$.  Two bounding disks are required to be disjoint.  
A regular neighborhood of the union of $G$ and the 
bounding disks is called an \emph{$s$--regular neighborhood} of $G$.  
\end{defi}
In particular, a $0$--special leaf with respect to $G$ is called a \emph{trivial leaf}.    
If a $Y_k$--graph $G$ in $M$ contains a $0$--special leaf $f$ with respect to $G$, 
then $M_G$ is diffeomorphic to $M$ \cite{habiro,GGP}.  

In the rest of the paper, a \emph{special leaf} is an $m$--special leaf for some unspecified integer $m$.\footnote{Note that in some 
literature \cite{GGP} the terminology `special leaf' is used to denote a $(-1)$--special leaf.  }   
The mention `with respect to' will be omitted when $G$ is clear from the context.  
\subsection{Statement of the result} \label{FFF}
Let $G$ be a $Y_n$--tree in a $3$--manifold $M$, $n\ge 2$.  It is well-known that, if 
$G$ contains a $(-1)$--special leaf, then 
  \begin{equation} \label{ZIP}
    M_G\sim_{Y_{n+1}} M. 
  \end{equation}
\noindent See \cite[Lemma E.21]{ohtsuki} for a proof for $M=S^3$, which can be generalized to our context.  
See also \cite[Lemma 4.9]{GGP}.    

We obtain the following generalization.  
\begin{thm} \label{table}
  Let $G$ be a $Y_n$--tree in a $3$--manifold $M$, with $n\ge 2$.  
  Let $l$ denote the number of special leaves with respect to $G$.  Then 
   \begin{enumerate}
      \item If $l<n$, then $M_G\sim_{Y_{n+l}} M.$ 
      \item If $l=n$, then $M_G\sim_{Y_{2n-1}} M.$ 
      \item If $l>n$, then $M_G\sim_{Y_{2n}} M.$    
   \end{enumerate}
\end{thm}
The proof is given in \fullref{prooftable}.  
In the next three subsections, we prove \fullref{table} in several important cases and provide 
a lemma which is used in \fullref{prooftable}.  
\subsection{The case of a tree with one special leaf} \label{l=1}
In this subsection, we prove \fullref{table} for $l=1$.  
\begin{lem} \label{onespecial}
  Let $G$ be a $Y_n$--tree in a $3$--manifold $M$, with $n\ge 2$. 
  Suppose that $G$ contains an $m$--special leaf ; $m\in \mathbf{Z}$. 
  Then $M_G\sim_{Y_{n+1}} M$.  
\end{lem}
\begin{proof}
 We first prove the lemma for all $m<0$, by induction.  
 As recalled in \fullref{FFF}, we already have the result for $m=-1$.  Now consider a $Y_n$--tree $G$ in $M$ with an 
 $m$--special leaf $f$, $m<0$.   
 Denote by $G'$ the clasper obtained by replacing $f$ by the union of a box $b$ and two edges $e_1$ and $e_2$ connecting $b$ 
 respectively to a $(-1)$--special leaf $f_1$ and a $(m+1)$--special leaf $f_2$ (both leaves being special with respect to $G'$).  
 By Habiro's move 7, $G'\sim G$.  
 Denote by $G_i$ the $Y_n$--tree in $M$ obtained from $G$ by replacing $f$ by $f_i$ ($i=1,2$).  
 By a zip construction, we have 
  $$ G'\sim (G_1\cup P), $$
 where $P$ satisfies $P\sim G_2$.  By \eqref{ZIP} it follows that $M_G\sim_{Y_{n+1}} M_{G_2}$.  
 The result then follows from the induction hypothesis.  
 
 Similarly, it would suffice to show the result for $m=1$ to obtain, by a similar induction, the result for all $m>0$.  
 For this, consider the case $m=0$.  In this case, $f$ is a trivial leaf and therefore $M_G\cong M$.  
 The same construction as above, with a $(-1)$--special leaf $f_1$ and a $1$--special leaf $f_2$,  
 shows that $M\sim_{Y_{n+1}} M_{G'}$, where $G'$ is a $Y_n$--tree in $M$ with a $1$--special leaf.  
 This concludes the proof.  
\end{proof}
\subsection{The case of a $Y_2$--tree} \label{n=2}
In this section, we prove \fullref{table} for $n=2$.  
The proof mainly relies on the following lemma.  
\begin{lem} \label{special}
   Let $G$ be a $Y_2$--tree in a $3$--manifold $M$ which contains two \mbox{$(-1)$--special} leaves which are connected to the same node.    
   Then $M_G\sim_{Y_{4}} M$.  
\end{lem}
\begin{proof}
 Denote by $w$ the node of $G$ which is connected to the two special leaves.  $w$ is connected by an edge to another node $v$.  
 By applying \fullref{move} at $v$, $G$ is equivalent, in an $s$--regular neighborhood, to a clasper $G'$ which is identical 
 to $G$, except in a $3$--ball where it is as depicted in \fullref{prooftwo} (a).  
 There, the node $w'$ corresponds to the node $w$ of $G$.  
 By \fullref{slide} and Habiro's move 6, we obtain the clasper depicted in \fullref{prooftwo} (b), 
 which is equivalent to the one depicted in \fullref{prooftwo} (c) by three applications of Habiro's move 12, 
 \fullref{lass} and an isotopy.  
 Denote by $G''$ this latter clasper.  As the figure shows, $G''$ contains a $Y_4$--subtree $T$.  
 Actually, $T$ is a `good input subtree' of $G''$, in the sense of 
\cite[Definition 3.13]{habiro}.     
    \begin{figure}[ht!]
    \bc
\labellist\small
\pinlabel $T$ [t] at 292 52
\hair6pt
\pinlabel (a) [t] at 32 13
\pinlabel (b) [t] at 132 13
\pinlabel (c) [t] at 278 13
\pinlabel $G'$ at 53 47
\pinlabel $G''$ at 209 47
\pinlabel $w'$ at 34 34
\endlabellist
    \includegraphics{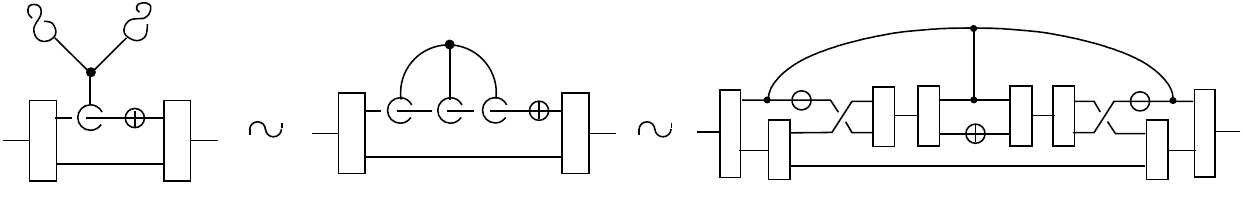} 
    \caption{} \label{prooftwo}
    \ec
    \end{figure}
 Denote by $\tilde{G''}$ the clasper obtained from $G''$ by inserting in each branch of $T$ 
 a pair of small Hopf-linked leaves.  By Habiro's move 2, $\tilde{G''}\sim G'$.  
 Denote by $\tilde{T}$ the $Y_4$--tree of $\tilde{G''}$ which corresponds to $T$.  
 By an application of the zip construction, we obtain 
 $M_{G''}\sim_{Y_{4}} M_{\tilde{G''}\setminus \tilde{T}}$.  
 Further, it follows from Habiro's moves 3 and 4 that $\tilde{G''}\setminus \tilde{T}\sim \emptyset$.   
\end{proof} 
The following technical lemma will allow us to generalize \fullref{special} to arbitrary special leaves.  
\begin{lem} \label{teck}
Let $G$ be a $Y_2$--tree in a $3$--manifold $M$ which contains two special leaves which are connected to the same node.  Then 
 $$ M_G\sim_{Y_{4}} M_{G_1\cup \tilde{G}_2}, $$  
where, for $i=1,2$, $G_i$ is obtained by replacing a $k$--special leaf of $G$ by a $k_i$--special leaf, such that $k_1+k_2=k$, and 
where $\tilde{G}_2$ is obtained from $G_2$ by an isotopy so that it is disjoint from $G_1$.  
\end{lem}
\begin{proof}
 Denote respectively by $f$ and $f'$ the $k$--special (resp.\ $k'$--special) leaf of $G$, $k,k'\in \mathbf{Z}$.  
 Just as in the proof of \fullref{onespecial}, we can use Habiro's moves 7 and the zip construction to see that 
 $G$ is equivalent, in an $s$--regular neighborhood, to the clasper $C_1$ of 
 \fullref{FT}, where $f_1$ is a $k_1$--special leaf and $f_2$ is a $k_2$--special leaf such that $k_1+k_2=k$.   
 Consider the leaf of $C_1$ denoted by $F$ in the figure.  
 By Habiro's move 12 at $F$, followed by two applications of Habiro's move 11, we have $C_1\sim C_2$, where $C_2$ is represented 
 in \fullref{FT}.  
    \begin{figure}[ht!]
    \bc
\begin{picture}(0,0)%
\includegraphics{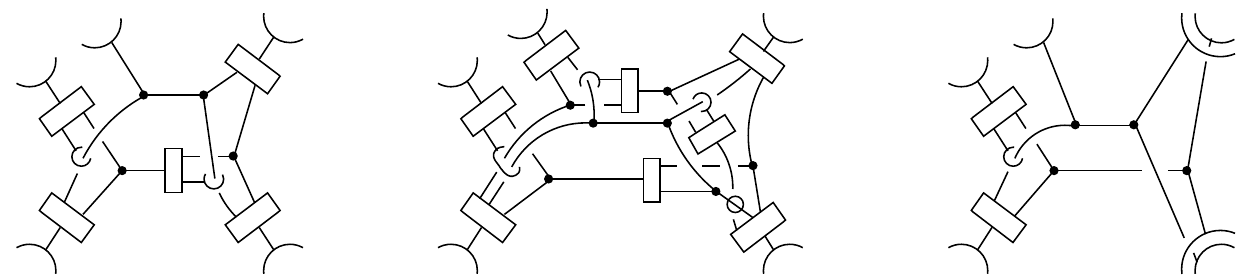}%
\end{picture}%
\setlength{\unitlength}{3947sp}%
\begingroup\makeatletter\ifx\SetFigFont\undefined%
\gdef\SetFigFont#1#2#3#4#5{%
  \reset@font\fontsize{#1}{#2pt}%
  \fontfamily{#3}\fontseries{#4}\fontshape{#5}%
  \selectfont}%
\fi\endgroup%
\begin{picture}(5920,1311)(34,-687)
\put(4820,489){\makebox(0,0)[lb]{\smash{\SetFigFont{8}{9.6}{\rmdefault}{\mddefault}{\updefault}
$f_2$
}}}
\put(2783,464){\makebox(0,0)[lb]{\smash{\SetFigFont{8}{9.6}{\rmdefault}{\mddefault}{\updefault}
$b$
}}}
\put(2374,546){\makebox(0,0)[lb]{\smash{\SetFigFont{8}{9.6}{\rmdefault}{\mddefault}{\updefault}
$f_2$
}}}
\put(4576,-361){\makebox(0,0)[lb]{\smash{\SetFigFont{8}{9.6}{\rmdefault}{\mddefault}{\updefault}
$b'$
}}}
\put(658,-662){\makebox(0,0)[lb]{\smash{\SetFigFont{8}{9.6}{\rmdefault}{\mddefault}{\updefault}
$C_1$
}}}
\put(2934,-662){\makebox(0,0)[lb]{\smash{\SetFigFont{8}{9.6}{\rmdefault}{\mddefault}{\updefault}
$C_2$
}}}
\put(5170,-662){\makebox(0,0)[lb]{\smash{\SetFigFont{8}{9.6}{\rmdefault}{\mddefault}{\updefault}
$C_4$
}}}
\put( 34,301){\makebox(0,0)[lb]{\smash{\SetFigFont{8}{9.6}{\rmdefault}{\mddefault}{\updefault}
$f_1$
}}}
\put( 56,-687){\makebox(0,0)[lb]{\smash{\SetFigFont{8}{9.6}{\rmdefault}{\mddefault}{\updefault}
$f'$
}}}
\put(4529,-687){\makebox(0,0)[lb]{\smash{\SetFigFont{8}{9.6}{\rmdefault}{\mddefault}{\updefault}
$f'$
}}}
\put(931,-377){\makebox(0,0)[lb]{\smash{\SetFigFont{8}{9.6}{\rmdefault}{\mddefault}{\updefault}
$F$
}}}
\put(347,489){\makebox(0,0)[lb]{\smash{\SetFigFont{8}{9.6}{\rmdefault}{\mddefault}{\updefault}
$f_2$
}}}
\put(4507,301){\makebox(0,0)[lb]{\smash{\SetFigFont{8}{9.6}{\rmdefault}{\mddefault}{\updefault}
$f_1$
}}}
\put(2057,301){\makebox(0,0)[lb]{\smash{\SetFigFont{8}{9.6}{\rmdefault}{\mddefault}{\updefault}
$f_1$
}}}
\put(2079,-687){\makebox(0,0)[lb]{\smash{\SetFigFont{8}{9.6}{\rmdefault}{\mddefault}{\updefault}
$f'$
}}}
\end{picture}
    \caption{} \label{FT}
    \ec
    \end{figure} 
 
 Consider the box $b$ of $C_2$ (see \fullref{FT}).  By applying Habiro's move 5 at $b$, $C_2$ is equivalent to a clasper containing 
 a $Y_3$--subtree $T$ and a $Y_1$--subtree $T'$ such that both $T$ and $T'$ contain a copy of $f_2$.  Denote by $C_3$ 
 the clasper obtained by replacing these two (linked) copies of $f_2$ by two $k_2$--special leaves.  
 By \fullref{crossingchange}, we have $M_{C_2}\sim_{Y_4} M_{C_3}$.  It follows from \fullref{onespecial} 
 and Habiro's move 5 that $M_{C_3}\sim_{Y_4} M_{C_4}$, where $C_4$ is as represented in \fullref{FT}.  
 By applying Habiro's move 5 at the box $b'$, $C_4$ is equivalent to a clasper containing 
 a $Y_2$--tree and a $Y_2$--subtree, each containing a copy of $f'$.  
 By \fullref{crossingchange}, $M_{C_4}\sim_{Y_4} M_{C_5}$, where $C_5$ is obtained by replacing these two (linked) copies of $f'$ 
 in $C_4$ by two $k'$--special leaves.  
 The result then follows from an isotopy and Habiro's move 3.  
 \end{proof}
We can now prove the case $n=2$ of \fullref{table}.  

Let $G$ be a $Y_2$--tree in a $3$--manifold $M$ with $l$ special leaves.  
If $l=0$, then the result is obvious. 
If $l=1$, \fullref{onespecial} implies that $M_G\sim_{Y_{3}} M$.  
If $l=2$, then $M_G\sim_{Y_{3}} M$ also follows from \fullref{onespecial}.  
It remains to prove the result when $l=3$ or $4$.  

Let $k,k'\in \mathbf{Z}$.  Denote by $G_{k,k'}$ a $Y_2$--tree in $M$ containing a $k$--special leaf $f$ 
and an $k'$--special leaf $f'$, both connected to the same node.  
Observe that it suffices to show that 
 \begin{equation} \label{EEE}
   M_{G_{k,k'}}\sim_{Y_{4}} M
 \end{equation} 
 If $k=k'=-1$, then \eqref{EEE} follows from \fullref{special}.  
Now, let us fix $k'=-1$.  Then we can show by induction that \eqref{EEE} holds for all $k<-1$.  
Indeed, consider some integer $m<-1$, and consider $G_{m,-1}$ in $M$.  
By \fullref{teck}, we have 
  $$ M_{G_{m,-1}}\sim_{Y_{4}} M_{C_1\cup C_2}, $$  
where $C_1$ contains two $(-1)$--special leaves connected to the same node, and where $C_2$ contains a $(-1)$--special leaf 
and an $m+1$--special leaf, both connected to the same node.  By \fullref{special} and the induction hypothesis, we 
thus obtain $M_{G_{m,-1}}\sim_{Y_{4}} M$. 

So we can now set $k'$ to be any negative integer, and prove \eqref{EEE} for all $k<-1$, by strictly the same induction.  

Similarly, it would suffice to show the result for $G_{1,1}$ to be able to prove \eqref{EEE} for all $k,k'>0$.  
Consider $G_{0,-1}$ in $M$.  In this case, $f$ is a trivial leaf and $M_{G_{0,1}}\cong M$.  
By applying \fullref{teck} at $f$, 
  $$ M_{G_{0,1}}\cong M\sim_{Y_{4}} M_{G_1\cup G_2}, $$ 
where $G_1$ (resp.\ $G_2$) contains a 
$(-1)$--special leaf and a $1$--special (resp.\ $(-1)$--special) leaf, both connected to the same node.  
It follows from \fullref{onespecial} that $M\sim_{Y_{4}} M_{G_1}$.  This proves \eqref{EEE} for $k=1$ and $k=-1$.  
We obtain \eqref{EEE} for $k=k=1$ similarly, by applying \fullref{teck} to $G_{0,1}$ in $M$.  
\subsection{The cutting lemma.}
Let $G$ be a $Y_n$--tree in $M$, with $n\ge 3$.  
By inserting a pair of small Hopf-linked leaves in an edge of $G$, we obtain a $Y_{n_1}$--tree $G_1$ and a $Y_{n_2}$--tree 
$G_2$ such that $n_1+n_2=n$ and $G_1\cup G_2\sim G$ (by Habiro's move 2).  See \fullref{figcut}.  
\begin{lem} \label{cut}
  Let $i=1,2$.   
  Suppose that, in a regular neighborhood $N_i$ of $G_i$, we have $(N_i)_{G_i}\sim_{Y_{k_i}} N_i$, with $k_1\ge 2$ and $k_2\ge 1$.  
  Then \be
   \item $M_G\sim_{Y_{k_1+2}} M$, if $G_2$ is a $Y_1$--tree containing at least one special leaf with respect to $G_1\cup G_2$,
\item $M_G\sim_{Y_{k_1+k_2}} M$, otherwise.\ee
\end{lem}
\begin{proof}
 Denote by $N$ an $s$--regular neighborhood of $G\sim G_1\cup G_2$.  
 Consider a $3$--ball $B$ in $M$ which intersects $N$ and $G_1\cup G_2$ as depicted in \fullref{figcut} (a).  
    \begin{figure}[ht!]
    \bc
\labellist\small
\pinlabel (a) [t] at 70 8
\pinlabel (b) [t] at 291 8
\pinlabel* $G_1$ [bl] at 27 57
\pinlabel* $G_2$ [br] at 112 57
\pinlabel* $N$ [tl] at 34 15
\pinlabel* $N'$ [tl] at 258 18
\pinlabel* $N''$ [tr] at 319 18
\pinlabel $G'_1$ <0pt,1pt> at 255 52
\pinlabel $G'_2$ at 325 52
\pinlabel $B$ [b] at 70 67
\pinlabel $B$ [b] at 290 67
\endlabellist
    \includegraphics{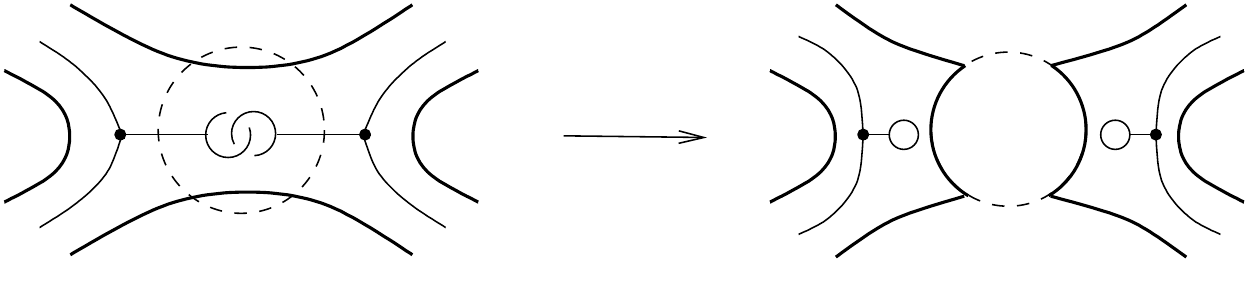} 
    \caption{} \label{figcut}
    \ec
    \end{figure}
 Denote by $N'$ and $N''$ the two connected components of the closure of $N\setminus (B\cap N)$.  
 By one crossing change and isotopy, we can homotop the two Hopf-linked leaves of $G_1\cup G_2$ into $N\setminus (B\cap N)$ so that,  
 if $G'_1\cup G'_2$ denotes the resulting clasper, we have $G'_1\subset N'$ and $G'_2\subset N''$.  See \fullref{figcut} (b).  
 Each of $G'_1$ and $G'_2$ contains a trivial leaf with respect to $G'_1\cup G'_2$, so we have 
 $G'_1\cup G'_2\sim \emptyset$ in $N$.  
 
 We now prove (1): suppose that $G_2$ contains one node and at least one special leaf with respect to $G_1\cup G_2$. 
 Denote by $f$ the leaf of $G_2$ which forms a Hopf link with a leaf of $G_1$.  
 By assumption, $G_1$ can be replaced by a $Y_{k_1}$--forest $F_1$ in an $s$--regular neighborhood $N_1$ so that $F_1\cup G_2\sim G$ 
 in $N$.  
 Consider a disk $d$ bounded by $f$ such that $d$ intersects transversally edges and leaves of components of $F_1$.  
 By a sequence of crossing changes, we can homotop these edges and leaves into $N'\subset N$: the clasper $G'$ obtained from 
 $F_1\cup G_2$ by this homotopy satisfies $G'\sim G'_1\cup G'_2\sim \emptyset$ in $N$.  
 So it would suffice to show that $M_{F_1\cup G_2}\sim_{Y_{k_1+2}} M_{G'}$.

 By \fullref{crossingchange}, we have $M_{F_1\cup G_2}\sim_{Y_{k_1+2}} M_{\tilde{F_1}\cup \tilde{G_2}}$, where 
 $\tilde{F_1}\cup \tilde{G_2}$ is obtained by `homotoping' into $N'$ all edges of $F_1$ and all $Y_{k}$--trees of $F_1$ with $k>k_1$.  
 Denote by $\tilde{f}$ the leaf of $\tilde{G_2}$ corresponding to $f$.  
 There is a sequence of crossing changes 
   $$ \tilde{F_1}\cup \tilde{G_2}=C_0\mapsto C_1\mapsto C_2\mapsto ...\mapsto C_{p-1}\mapsto C_p=G', $$
 where, for each $1\le k\le p$, $C_k$ is obtained from $C_{k-1}$ by one crossing change between $\tilde{f}$ and a leaf $l$ of a 
 $Y_{k_1}$--tree $T_k$ of $\tilde{F_1}$.\footnote{
 Here, abusing notations, we still denote by $\tilde{f}$, $\tilde{G_2}$ and $\tilde{F_1}$ the corresponding elements in $C_k$, 
 for all $k\ge 1$.  }
 By \fullref{crossingchange}, we have $M_{C_{k}}\sim_{Y_{k_1+2}} M_{C_{k-1}\cup H_k}$, where $H_k$ is a $Y_{k_1+1}$--tree obtained 
 by connecting the edges of $\tilde{G_2}$ and $T_k$ attached to $\tilde{f}$ and $l$ respectively.  In particular,   
 $H_k$ contains a special leaf with respect to $C_{k-1}\cup H_k$.  So by \fullref{onespecial}, we have 
 $M_{C_{k}} \sim_{Y_{k_1+2}} M_{C_{k-1}}$.  
 It follows that $M_{\tilde{F_1}\cup \tilde{G_2}} \sim_{Y_{k_1+2}} M_{G'}$, which concludes the proof of (1).  
 
 The proof of (2) is simpler, and left to the reader.  It uses exactly the same arguments as above, 
 by considering the $Y_{k_i}$--forest $F_i$ ($i=1,2$) in an $s$--regular neighborhood $N_i$ of $G_i$ such that  
 $F_1\cup F_2\sim G$ in $N$.  
\end{proof}
\subsection[Proof of \ref{table}]{Proof of \fullref{table}} \label{prooftable}
Suppose that $G$ is a $Y_n$--tree in $M$ with $l$ special leaves ; $n\ge 2$, $l\ge 0$.  
\subsubsection{The case $l<n$}  \label{l<n}
In this case, it is necessary to reduce the problem to linear trees.  
We have the following. 
\begin{claim} \label{cl}
Let $1\le p\le l$ be an integer. 
Pick two non-special leaves $f_1$ and $f_2$ of $G$.  
Then we have, by successive applications of the IHX relation, 
 $$ M_G\sim_{Y_{n+p}} M_{L_p},  $$ 
where $L_p$ is a union of disjoint linear $Y_k$--trees with $n\le k\le n+p-1$ such that 
 \bi    
   \item the ends of each linear tree are parallel copies of $f_1$ and $f_2$, 
   \item each $Y_k$--tree contains $(n+l-k)$ special leaves with respect to $L_p$.  
 \ei 
\end{claim}
\begin{proof}[Proof of the claim]
 The claim is proved by induction on $p$.  
 Observe that we can use the IHX relation to replace $T$ by a union $L_1$ 
 of linear $Y_k$--trees whose ends are parallel copies of $f_1$ and $f_2$.  \fullref{crossingchange} (1) ensures that each tree has $l$ special leaves 
 with respect to $L_1$.  This proves the case $p=1$.
 Now assume the claim for some $p\ge1$: $M_T\sim_{Y_{n+p}} M_{L_p}$, where 
 $L_p$ is as described above.  By assumption, this equivalence comes from 
 \fullref{ihx}, so we can apply \fullref{corihx}.  There exists a 
 union $F$ of disjoint (possibly non linear) $Y_{n+p}$--trees such that
 $M_T\sim_{Y_{n+p+1}} M_{L_p\cup F}$.  For each tree $T$ in $F$, its ($n+p+2$) 
 leaves are obtained by taking the leaves of a $Y_{n+p-1}$--tree in $L_p$ and 
 adding a parallel copy of one of them.  If this additional leaf is a copy of 
 a special leaf $f$ (with respect to $L_p$), the two (linked) copies of $f$ in 
 $T$ are not special leaves with respect to $L_p\cup F$.  This shows that each 
 tree in $F$ contains at least ($l-p$) special leaves with respect to $L_p\cup F$.  
 Note that each such tree also contains (at least) a copy of $f_1$ and $f_2$.  
 So by \fullref{ihx} we have $M_{L_p\cup F}\sim_{Y_{n+p+1}} M_{L_{p+1}}$, where 
 $L_{p+1}$ is of the desired form.  
\end{proof}
It follows from Claim \ref{cl} that 
 \begin{equation*} 
   M_T\sim_{Y_{n+l}} M_L, 
 \end{equation*}
where $L$ is a union of linear $Y_k$--trees with $n\le k\le n+l-1$,    
each such linear $Y_k$--tree containing (at least) $(n+l-k)$ special leaves with respect to $L$, and whose ends are non-special leaves. 

So it suffices to prove the case $l<n$ of \fullref{table} for linear $Y_n$--trees whose ends are non-special leaves. 
We proceed by induction on $n$.  

For $n=2$, the statement follows from \fullref{n=2}.  

Now, assume that the statement holds true for all $k<n$, and consider a linear $Y_n$--tree $G$ whose ends are two non-special leaves.  
Insert a pair of small Hopf-linked leaves in an edge of $G$ such that it produces a union of two linear trees $G_1\cup G_2\sim G$ 
with $degG_1=n_1$ and $degG_2=n_2$.  Denote respectively by $l_1$ and $l_2$ the number of special leaves with respect to $G_1\cup G_2$ 
in $G_1$ and $G_2$.  We have $n_1+n_2=n$ and $l_1+l_2=l$.  Denote also by $N_1$ an $s$--regular neighborhood of $G_1$.  
 \begin{itemize}
  \item If we can choose $n_2=1$ and $l_2=1$, then $n_1=n-1$ and $l_1=l-1$. So $l_1<n_1$ and by the induction 
        hypothesis we have $(N_1)_{G_1}\sim_{Y_{n+l-2}} N_1$ ($G_1$ is indeed linear).   As $G_2$ contains one special 
        leaf with respect to $G_1\cup G_2$, we obtain the result by \fullref{cut} (1).  
  \item Otherwise, then $l<n-1$, and we can choose $G_2$ such that $n_2=1$ and $l_2=0$ (that is, $G_2$ contains one node connected to 
        $2$ non-special leaves).  
        As $l_1=l< n_1=n-1$, we have $(N_1)_{G_1}\sim_{Y_{n+l-1}} N_1$ (by the induction hypothesis), and the 
        result follows from \fullref{cut} (2).  
 \end{itemize}    
This completes the proof of the case $l<n$.  
\subsubsection{The case $l\ge n$}  
The case $l=n$ follows immediately from the case $l=n-1$, by regarding one of the special leaves as a leaf.  

We prove the case $l=n+1$ by induction on the degree $n$.  
The case $n=2$ was proved in \fullref{n=2}.  
Consider a $Y_n$--tree $G$ with $l\ge n$ special leaves.  As in \fullref{l<n}, insert 
a pair of Hopf-linked leaves in an edge of $G$ so that we obtain a union of two trees $G_1\cup G_2\sim G$ 
with $degG_1=n-1$ and $degG_2=1$.  
Denote respectively by $l_1$ and $l_2$ the number of special leaves with respect to $G_1\cup G_2$ in $G_1$ and $G_2$.  
There are two cases, depending on whether $l_2=1$ or $2$.  
 \bi
  \item If $l_2=1$, then $l_1=n=n_1+1$, and thus, by the induction hypothesis we have 
        $(N_1)_{G_1}\sim_{Y_{2n-3}} N_1$ in an $s$--regular neighborhood $N_1$ of $G_1$.  
        The result follows from \fullref{cut} (1).  
  \item If $l_2=2$, then $l_1=n-1=n_1$.  It thus follows from the case $l=n$ of \fullref{table} 
        that $(N_1)_{G_1}\sim_{Y_{2n-3}} N_1$ in an $s$--regular neighborhood $N_1$ of $G_1$.  
        The result then follows as above from \fullref{cut} (1).  
 \ei
 
The case $l=n+2$ follows from the case $l=n+1$ by regarding one of the special leaves as a leaf.  
\subsection[Some special cases for \ref{table}]{Some special cases for \fullref{table}} \label{mariethomas}
We have the following improvement of \fullref{table} for linear trees having only $(-1)$--special leaves.  
\begin{prop} \label{linear2n+1}
  Let $G$ be a linear $Y_n$--tree in a $3$--manifold $M$, $n\ge 2$, such that all its leaves are $(-1)$--special leaves.  
  Then in an $s$--regular neighborhood $N$ of $G$ (which is a $3$--ball in $M$) we have 
    $$ N_G\sim_{Y_{2n+1}} N_{\Theta_n}, $$ 
  where $\Theta_n$ is the connected $Y_{2n}$--graph without leaves depicted in \fullref{thetan}
   \begin{figure}[ht!]
   \bc
    \includegraphics{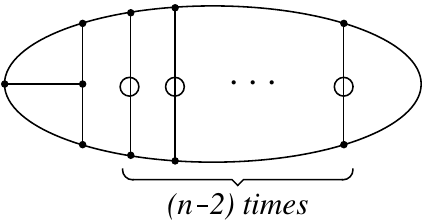}
    \caption{The $Y_{2n}$--graph $\Theta_n$} \label{thetan}
   \ec
   \end{figure}
\end{prop}
\begin{rem}
Note that ``$\sim_{Y_{2n+1}}$'' in \fullref{linear2n+1} can be replaced by ``$\sim_{Y_{2n+2}}$''.  
This follows from the fact that if two integral homology balls are $Y_{2n+1}$--equivalent then they
are $Y_{2n+2}$--equivalent ($n>1$).
\end{rem}
The proof of \fullref{linear2n+1} uses rather involved calculus of claspers, and is therefore postponed to 
\fullref{zeproof}.  
Note that this result is not needed for the rest the paper.  
A reader who is not too comfortable with claspers (but who nevertheless reached this point) may thus safely skip this proof.  

Also, one can check that if $G$ is a $Y_n$--tree in a $3$--manifold $M$ with $n$ special leaves, we have 
  \begin{equation} \label{opla} 
   M_{G}\sim_{Y_{2n}} M 
  \end{equation}
in the two following situations: 
\bi
\item $G$ contains a $2k$--special leaf, for some integer $k$.  
\item The homology class in $H_1(M;\mathbf{Z} / 2\mathbf{Z})$ of a non-special leaf of $G$ is zero.\footnote{  
      This fact was pointed out to the author by Kazuo Habiro.  }
      In particular, \eqref{opla} always holds if $M=S^3$.  
\ei 
\section{$Y_k$--equivalence for $3$--manifolds obtained by surgery along Brunnian links} \label{+1surgery}
In this section, we prove Theorems \ref{2n-2} and \ref{2n-1}.  
The proofs use a characterization of Brunnian links in terms of claspers due to Habiro, and independently 
to Miyazawa and Yasuhara, which involves the notion of $C^a_k$--equivalence.  
Let us first recall from \cite{hbrun} the definition and some properties of this equivalence relation.  
\subsection{$C^a_k$--equivalence} \label{sec:cak-equivalence}
\begin{defi} \label{r9}
  Let $L$ be an $m$--component link in a $3$--manifold $M$.  For $k\ge m-1$,
  a {\em $C^a_k$--tree} for $L$ in $M$ is a $C_k$--tree $T$ for $L$ in
  $M$, such that
  \begin{enumerate}
  \item all the strands intersecting a given disk-leaf of $T$ are from the same component of $L$, 
  \item $T$ intersects \emph{all} the components of $L$.
  \end{enumerate}
\end{defi}

A (simple) \emph{$C^a_k$--forest} $L$ is a clasper consisting only of (simple) $C^a_k$--tree for $L$.  

A {\em $C^a_k$--move} on a link is surgery along a $C^a_k$--tree.  The
{\em $C^a_k$--equivalence} is the equivalence relation on links
generated by $C^a_k$--moves.  

The main tool in the proofs of Theorems \ref{2n-2} and \ref{2n-1} is the following.  
\begin{thm}{\rm\cite{hbrun,MY}}\label{brunian}\qua
Let $L$ be an $(n+1)$--component link in $S^3$.  $L$ is Brunnian if and only if it is 
\emph{$C^a_n$--equivalent} to the $(n+1)$--component unlink $U$.
\end{thm}

In the proof of \fullref{2n-1}, we will also need the next result.  
\begin{thm}[\cite{MY}, see also \cite{hm1}] \label{lh}
Two $(n+1)$--component Brunnian links in $S^3$ are link-homotopic 
if and only if they are $C^a_{n+1}$--equivalent.  
\end{thm}
\noindent Note that this statement does not appear explicitly in 
\cite{MY}.  However, it is implied by the proof of \cite[Theorem 3]{MY}.  An alternative proof was given subsequently by Habiro 
and the author \cite{hm1}.  
\subsection[Proof of \ref{2n-2}]{Proof of \fullref{2n-2}} \label{proof2n-2}
Let $m=(m_1,...,m_{n+1})\in \mathbf{Z}^{n+1}$, $n\ge 2$ and let $L$ be an $(n+1)$--component Brunnian link in a $3$--manifold $M$.   
By \fullref{brunian}, $L$ is $C^a_n$--equivalent to an $(n+1)$--component unlink $U$ in $M$.  So by \cite[Lemma 7]{hbrun} 
there exists a simple $C^a_n$--forest $F=T_1\cup ...\cup T_p$ for $U$ such that $L\cong U_F$.  
We thus have
  $$ M_{(L,m)}\cong M_{G_m(F)}, $$ 
where $G_m(F)$ is the clasper obtained from $F$ by performing $\frac{1}{m_i}$--framed surgery along the $i^{th}$ component $U_i$ 
of $U$ for all $1\le i\le n+1$.  
Indeed, $\frac{1}{m_i}$--surgery along an unknot does not change the diffeomorphism type of $M$, and can be regarded as a move on 
claspers in $M$.  
Observe that  
$\frac{1}{m_i}$--surgery along $U_i$ turns each disk-leaf of $F$ intersecting $U_i$ into a $(-m_i)$--framed unknot (here, we forget 
the bounding disk).  Thus $\frac{1}{m}$--surgery along $U$ turns each $C^a_n$--tree $T_j$ of $F$ into a $Y_{n-1}$--tree $G_j$ in $M$.  
However, the $(n+1)$ corresponding leaves of $G_j$ might not be special leaves with respect to $G_m(F)$, as they can be linked 
with the leaves of other components of $G_m(F)$.  \fullref{crossingchange} (1) can be used to unlink these leaves 
`up to $Y_{2n-2}$--equivalence'.  Namely, \fullref{crossingchange} implies that $M_{G_m(F)}\sim_{Y_{2n-2}} M_{\tilde{G}_m(F)}$, 
where $\tilde{G}_m(F)$ is a union of $Y_{n-1}$--trees, each containing $(n+1)$ special leaves with respect to $\tilde{G}_m(F)$.  
The result then follows from \fullref{table}.  
\subsection[Proof of \ref{2n-1}]{Proof of \fullref{2n-1}} \label{proof2n-1}
Let $L$ and $L'$ be two link-homotopic $(n+1)$--component Brunnian links in $M$, and let $U$ denote an $(n+1)$--component unlink $U$ 
in $M$.  
By Theorems \ref{brunian} and \ref{lh}, and \cite[Lemma 7]{hbrun}, there exists a simple $C^a_{n+1}$--forest $F=T_1\cup ...\cup T_p$ 
and a simple $C^a_{n}$--forest $F'=T'_1\cup ...\cup T'_q$ for $U$ such that $L'\cong U_{F'}$ and $L\cong U_{F\cup F'}$.  

For all $j$, denote by $G'_j$ (resp.\ $G_j$) the $Y_{n-1}$--tree (resp.\ $Y_n$--tree) obtained from $T'_j$ (resp.\ $T_j$) 
by $\frac{1}{m}$--surgery along $U$.  
By \fullref{crossingchange}, 
  $$ M_{(L,+1)}\sim_{Y_{2n-1}} M_{G'_1\cup ...\cup G'_q}\sharp S^3_{G_1}\sharp ... \sharp S^3_{G_p}\cong 
     M_{(L',+1)}\sharp S^3_{G_1}\sharp ... \sharp S^3_{G_p}. $$
So proving that $S^3_{G_i}\sim_{Y_{2n-1}} S^3$ for all $1\le i \le p$ would imply the theorem.  

By strictly the same arguments as in \fullref{proof2n-2}, the $Y_{n}$--tree $G_i$ contains at least  
$n$ special leaves, for all $1\le i\le p$.  
So \fullref{table} implies that $S^3_{G_i}\sim_{Y_{2n-1}} S^3$.  
\section{Trivalent diagrams and Goussarov--Vassiliev invariants for Brunnian links} \label{GV}
In this section, we recall some results proved by Habiro and the author in a previous paper \cite{hm2}.  These, together with the 
two theorems shown in \fullref{+1surgery}, will allow us to prove \fullref{homspheres} in the next section. 
\subsection{Trivalent diagrams} \label{trivalent}
A \emph{trivalent diagram} is a finite graph with trivalent vertices, each vertex being equipped with a cyclic order on the 
three incident edges.  The \emph{degree} of a trivalent diagram is half the number of vertices.  

\begin{figure}[ht!]
  \bc
\begin{picture}(0,0)%
\includegraphics[scale=1.2]{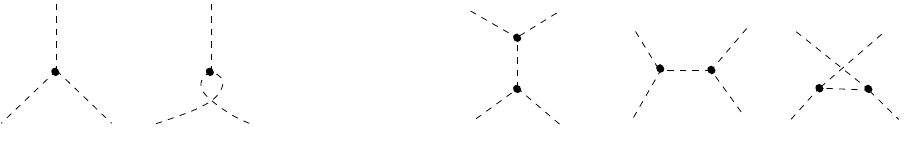}%
\end{picture}%
\setlength{\unitlength}{2367sp}%
\def\SetFigFont#1#2#3#4#5{\small}%
\begin{picture}(8734,1544)(2121,-6159)
\put(9353,-5238){\makebox(0,0)[lb]{\smash{\SetFigFont{6}{7.2}{\rmdefault}{\mddefault}{\updefault}{\color[rgb]{0,0,0}$+$}%
}}}
\put(10855,-5223){\makebox(0,0)[lb]{\smash{\SetFigFont{6}{7.2}{\rmdefault}{\mddefault}{\updefault}{\color[rgb]{0,0,0}$=0$}%
}}}
\put(7645,-5238){\makebox(0,0)[lb]{\smash{\SetFigFont{6}{7.2}{\rmdefault}{\mddefault}{\updefault}{\color[rgb]{0,0,0}$-$}%
}}}
\put(3212,-5238){\makebox(0,0)[lb]{\smash{\SetFigFont{6}{7.2}{\rmdefault}{\mddefault}{\updefault}{\color[rgb]{0,0,0}$+$}%
}}}
\put(4624,-5223){\makebox(0,0)[lb]{\smash{\SetFigFont{6}{7.2}{\rmdefault}{\mddefault}{\updefault}{\color[rgb]{0,0,0}$=0$}%
}}}
\put(8647,-6128){\makebox(0,0)[lb]{\smash{\SetFigFont{6}{7.2}{\rmdefault}{\mddefault}{\updefault}{\color[rgb]{0,0,0}IHX}%
}}}
\put(3242,-6159){\makebox(0,0)[lb]{\smash{\SetFigFont{6}{7.2}{\rmdefault}{\mddefault}{\updefault}{\color[rgb]{0,0,0}AS}%
}}}
\end{picture}
    \caption{The AS and IHX relations} \label{ASIHX}
  \ec
  \end{figure}
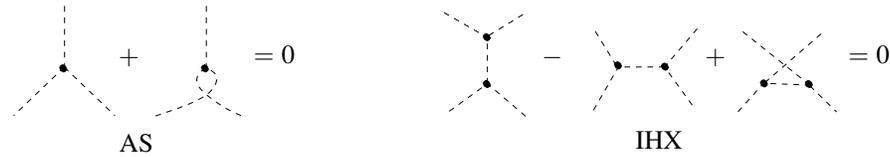

For $k\ge 0$, let $\mathcal{A}_k(\emptyset)$ denote the $\mathbf{Z}$--module generated by trivalent diagrams of degree 
$k$, subject to the \emph{AS and IHX relations}, see \fullref{ASIHX}.

Denote by $\mathcal{A}^c_k(\emptyset)$ the $\mathbf{Z}$--submodule of $\mathcal{A}_k(\emptyset)$ generated by \emph{connected} 
trivalent diagrams. 
\subsection{The Brunnian part of the Goussarov--Vassiliev filtration} \label{brunnianpart}
Denote by $\mathbf{Z}\mathcal{L}(n)$ the free $\mathbf{Z}$--module generated by the set of isotopy classes of $n$--component links 
in $S^3$, and denote by $J_k(n)$ the $\mathbf{Z}$--submodule of $\mathbf{Z}\mathcal{L}(n)$ generated by elements of the form 
  $$ [L;C_1,...,C_p]:=\sum_{S\subseteq \{ C_1, ... ,C_p\} } (-1)^{|S|} L_S, $$
where $L$ is an $n$--component link in $S^3$, and where the $C_i$ ($1\le i\le p$) are disjoint $C_{k_i}$--trees for $L$ such that 
$k_1+...+k_p=k$.  The sum runs over all the subsets $S$ of $\{ C_1, ... ,C_p\}$ and $|S|$ denotes the cardinality 
of $S$.  
The descending filtration 
  $$ \mathbf{Z}\mathcal{L}(n)=J_0(n)\supset J_1(n)\supset J_2(n)\supset ... $$
coincides with the \emph{Goussarov--Vassiliev filtration} \cite{habiro}.  
  
Denote by $\overline{J}_k(n)$ the graded quotient $J_k(n) / J_{k+1}(n)$.  
\begin{defi}
The \emph{Brunnian part} $\Br(\overline{J}_{2n}(n+1))$ of the $2n^{th}$ graded quotient $\overline{J}_{2n}(n+1)$ is 
the $\mathbf{Z}$--submodule generated by elements $[L-U]_{J_{2n+1}}$ where $L$ is an $(n+1)$--component Brunnian link. 
\end{defi}

\begin{figure}[ht!]
  \bc
\begin{picture}(0,0)%
\includegraphics{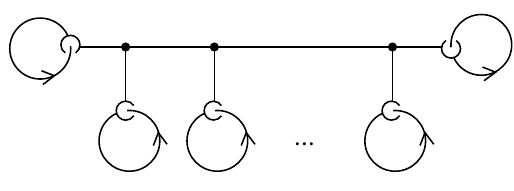}%
\end{picture}%
\setlength{\unitlength}{3947sp}%
\begingroup\makeatletter\ifx\SetFigFont\undefined%
\gdef\SetFigFont#1#2#3#4#5{%
  \reset@font\fontsize{#1}{#2pt}%
  \fontfamily{#3}\fontseries{#4}\fontshape{#5}%
  \selectfont}%
\fi\endgroup%
\begin{picture}(2465,915)(1,-406)
\put(1201,389){\makebox(0,0)[lb]{\smash{\SetFigFont{10}{12.0}{\rmdefault}{\mddefault}{\updefault}{\color[rgb]{0,0,0}$T_{\sigma}$}%
}}}
\put(226,-361){\makebox(0,0)[lb]{\smash{\SetFigFont{10}{12.0}{\rmdefault}{\mddefault}{\updefault}{\color[rgb]{0,0,0}$U_{\sigma(1)}$}%
}}}
\put(2026,-361){\makebox(0,0)[lb]{\smash{\SetFigFont{10}{12.0}{\rmdefault}{\mddefault}{\updefault}{\color[rgb]{0,0,0}$U_{\sigma(n-1)}$}%
}}}
\put(1201,-361){\makebox(0,0)[lb]{\smash{\SetFigFont{10}{12.0}{\rmdefault}{\mddefault}{\updefault}{\color[rgb]{0,0,0}$U_{\sigma(2)}$}%
}}}
\put(  1, 14){\makebox(0,0)[lb]{\smash{\SetFigFont{10}{12.0}{\rmdefault}{\mddefault}{\updefault}{\color[rgb]{0,0,0}$U_{n+1}$}%
}}}
\put(2176, 14){\makebox(0,0)[lb]{\smash{\SetFigFont{10}{12.0}{\rmdefault}{\mddefault}{\updefault}{\color[rgb]{0,0,0}$U_n$}%
}}}
\end{picture}
    \caption{The simple linear $C^a_{n}$--tree $T_{\sigma}$} \label{Tsigma}
  \ec
  \end{figure}
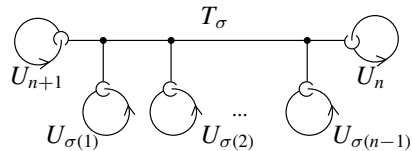

As outlined in \cite[Section 7]{hm1}, $\Br(\overline{J}_{2n}(n+1))$ is spanned over $\mathbf{Z}$ 
by elements 
  $$ \frac{1}{2}[U;T_{\sigma}\cup \tilde{T}_{\sigma}]\textrm{ and } 
     [U;T_{\sigma}\cup \tilde{T}_{\sigma'}],\quad  \textrm{ for $\sigma\ne \sigma' \in S_{n-1}$,}$$  
where, for all $\sigma,\sigma'$ in the symmetric group $S_{n-1}$, $T_{\sigma}$ is the simple linear $C^a_{n}$--tree for the 
$(n+1)$--component unlink $U$ depicted in \fullref{Tsigma}, and $\tilde{T}_{\sigma'}$ is obtained from $T_{\sigma'}$ by a small 
isotopy so that it is disjoint from $T_{\sigma}$.  (Here $\frac{1}{2}[U;T_{\sigma}\cup \tilde{T}_{\sigma}]$ means an element 
$x\in\Br(\overline{J}_{2n}(n+1))$such that $2x=[U;T_{\sigma}\cup \tilde{T}_{\sigma}]$.  Existence of such an
element is shown in \cite{hm1}.)
\subsection{The map $h_n\co  \mathcal{A}^c_{n-1}(\emptyset)\rightarrow \Br(\overline{J}_{2n}(n+1))$} 
\label{structure}
Connected trivalent diagrams allow us to describe the structure of $\Br(\overline{J}_{2n}(n+1))$.  
For $n\ge 2$, we have a map 
  $$  h_n\co  \mathcal{A}^c_{n-1}(\emptyset) \longrightarrow \overline{J}_{2n}(n+1) $$
defined as follows.  Given a connected trivalent diagram $\Gamma\in \mathcal{A}^c_{n-1}(\emptyset)$, 
insert $n+1$ ordered copies of $S^1$ in the edges of $\Gamma$, in an arbitrary way.  
The result is a strict unitrivalent graphs $D_{\Gamma}$ of degree $2n$ on the disjoint union of 
$(n+1)$ copies of $S^1$ (see \cite{BN}).  Next, `realize' this unitrivalent graph by a graph clasper.  
Namely, replace each univalent vertex (resp.\ trivalent vertex, edge) of $D_{\Gamma}$ with a disk-leaf 
(resp.\ node, edge), these various subsurfaces being connected as prescribed by the graph $D_{\Gamma}$.     
Denote by $C(D_{\Gamma})$ the resulting graph clasper for the $(n+1)$--component unlink $U$.  
Then 
  $$ h_n(\Gamma):=[U-U_{C(D_{\Gamma})}]_{J_{2n+1}}\in \overline{J}_{2n}(n+1). $$  
For $n\ge 3$, the image of $h_n$ is the Brunnian part $\Br(\overline{J}_{2n}(n+1))$ of $\overline{J}_{2n}(n+1)$, and 
  $$  h_n\oq\co  \mathcal{A}^c_{n-1}(\emptyset)\oq \longrightarrow \Br(\overline{J}_{2n}(n+1))\oq $$
is an isomorphism.    
\section{Finite type invariants of integral homology spheres} \label{FTIZHS}
\subsection{The Ohtsuki filtration for integral homology spheres} \label{sec:ohtsuki} 
Let $\mathcal{M}$ denote the free $\mathbf{Z}$--module generated by the set of orientation-preserving homeomorphism classes of 
integral homology spheres. The definition of the Ohtsuki filtration uses algebraically split, unit-framed links.  
For the purpose of the present paper, it is however more convenient to use a definition using claspers, due to Goussarov and 
Habiro \cite{GGP,goussarov,habiro}.  
For $k\ge 0$, let $\mathcal{M}_k$ denote the $\mathbf{Z}$--submodule of $\mathcal{M}$ generated by elements of the form
  $$ [M;G_1,...,G_p]:=\sum_{S\subseteq \{ G_1,...,G_p\} } (-1)^{|S|} M_S, $$
where $M$ is an integral homology sphere, and where the $G_i$ ($1\le i\le p$) are disjoint $Y_{k_i}$--graphs in $M$ such that 
$k_1+...+k_p=k$.  The sum runs over all the subsets $S$ of $\{ G_1,...,G_p\}$ and $|S|$ denotes the cardinality 
of $S$. 

The descending filtration of $\mathbf{Z}$--submodules 
  $$ \mathcal{M}=\mathcal{M}_0\supset \mathcal{M}_1\supset \mathcal{M}_2\supset ... $$ 
is equal to the Ohtsuki filtration after re-indexing and tensoring by 
$\mathbf{Z}[1/2]$ \cite{GGP,goussarov,habiro}.  
Another alternative definition was previously given by Garoufalidis and Levine using `blinks' \cite{GL}.  
%
\subsection{The connected part of the Ohtsuki filtration} \label{varphi}
\begin{figure}[ht!]
  \bc
    \includegraphics{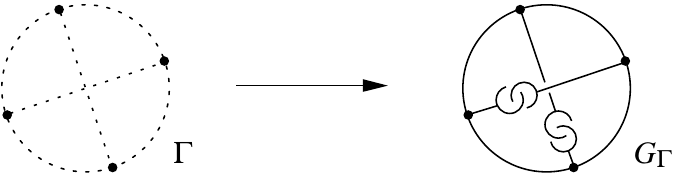}
    \caption{Realizing a trivalent diagram in $S^3$} \label{realizing}
  \ec
  \end{figure}

Let $\overline{\mathcal{M}}_{2k}$ denote the graded quotient $\mathcal{M}_{2k} / \mathcal{M}_{2k+1}$.  
\begin{defi}
The \emph{connected} part $\Co(\overline{\mathcal{M}}_{2k})$ of $\overline{\mathcal{M}}_{2k}$ 
is the $\mathbf{Z}$--submodule of $\overline{\mathcal{M}}_{2k}$ generated by elements $[S^3;G]_{\mathcal{M}_{2k+1}}$ 
where $G$ is a $Y_{2k}$--graph (in particular, $G$ is \emph{connected}).  
\end{defi}
For $k\ge 1$, there is a well-defined \emph{surgery map} 
  $$ \varphi_k\co  \mathcal{A}_k(\emptyset) \longrightarrow \overline{\mathcal{M}}_{2k}, $$
which maps each trivalent diagram $\Gamma=\Gamma_1\cup ...\cup \Gamma_p$ to $[S^3;G_{\Gamma_1},...,G_{\Gamma_p}]$, where $G_{\Gamma_i}$ 
is a connected clasper obtained by `realizing' the diagram $\Gamma_i$ in $S^3$ as depicted in \fullref{realizing}.  
The image $\varphi_k(\Gamma)$ of a degree $k$ trivalent diagram $\Gamma$ in $\overline{\mathcal{M}}_{2k}$ by $\varphi_k$ does not 
depend on the embeddings $G_{\Gamma_i}$ in $S^3$ (\cite{habiro}, see also \cite[page 320]{ohtsuki}).  
Note that $\varphi_k$ is a reconstruction, using claspers, of a map defined previously by Garoufalidis and Ohtsuki \cite{GO}.  
The homomorphism $\varphi_k\oz$ is surjective, and it is an isomorphism when tensoring by $\mathbf{Q}$, with inverse given by the 
LMO invariant \cite{Le}.  

It can be easily checked using the arguments of \cite{GGP} that 
$\varphi_k(\mathcal{A}^c_k(\emptyset))=\Co(\overline{\mathcal{M}}_{2k})$.  
We thus have an isomorphism 
  $$ \varphi_k\oq\co  \mathcal{A}^c_k(\emptyset)\oq \overset{\simeq}\longrightarrow \Co(\overline{\mathcal{M}}_{2k})\oq $$
induced by the surgery map $\varphi_k$.  
\subsection[The map alpha_k]{The map $\alpha_k\co  \Co(\overline{\mathcal{M}}_{2k})\longrightarrow \overline{\mathcal{S}}_{2k} $}  \label{sec:Sk}
Let $\mathcal{S}_k$ denote the set of integral homology spheres which are $Y_k$--equivalent 
to $S^3$, and denote by $\overline{\mathcal{S}}_{k}$ the quotient $\mathcal{S}_k / \sim_{Y_{k+1}}$.  
The connected sum induces an abelian group structure on $\overline{\mathcal{S}}_{k}$.  

As recalled in the introduction, $\overline{\mathcal{S}}_{2k+1}=0$ for all $k\ge 1$. 
$\overline{\mathcal{S}}_{2k}$ is generated by the elements $S^3_G$, where $G$ is a $Y_{2k}$--graph 
in $S^3$ (for $k=0$, we have $\overline{\mathcal{S}}_{1}=\mathbf{Z} / 2\mathbf{Z}$).  
There is a surjective homomorphism of abelian groups 
  $$ \phi_k\co  \mathcal{A}^c_{k}(\emptyset) \longrightarrow \overline{\mathcal{S}}_{2k} $$
defined by $\phi_k(\Gamma):=[S^3_{G_{\Gamma}}]_{Y_{2k+1}}$, 
where $G_{\Gamma}$ is a topological realization of the diagram $\Gamma$ as in the 
definition of $\varphi_k$ (see \fullref{realizing}).  
It is well known that $\phi_k$ is well-defined (see the proof of \cite[Theorem E.20]{ohtsuki}).  

The map $\phi_k$ is an isomorphism over the rationals.  This is shown by using the primitive part of the LMO invariant $z^{LMO}$
\cite[pages 329--330]{ohtsuki}.  

Let 
  $$ \alpha_k\co  \Co(\overline{\mathcal{M}}_{2k})\longrightarrow \overline{\mathcal{S}}_{2k} $$
be the map defined by 
  $$ \alpha_k([S^3;G]_{\mathcal{M}_{2k+1}})=[S^3_G]_{Y_{2k+1}}. $$ 
\noindent The fact that $\alpha_k$ is well-defined follows from standard arguments of clasper theory, 
and is well known to experts.   

The following is clear from the above definitions.  
  \begin{lem} \label{e2}
    The following diagram commutes for all $k\ge 1$:
     \begin{equation*}
      \xymatrix{
        \mathcal{A}^c_{k}(\emptyset)
        \ar[d]_{\varphi_k}
        \ar[dr]^{\phi_k}\\
        \Co(\overline{\mathcal{M}}_{2k})
        \ar[r]_{\quad \alpha_k\quad }&\overline{\mathcal{S}}_{2k}
      }
    \end{equation*}
  \end{lem}
As a consequence, $\alpha_k$ is an isomorphism over the rationals.  
\subsection{The map $\lambda_n$} \label{lambda}
For simplicity, we work over the rationals in the rest of this section.  

Let $n\ge 2$.  Denote by $B_{n+1}$ the set of isotopy classes of $(n+1)$--component Brunnian links in $S^3$.  
Define a linear map 
  $$ \tilde{\lambda}_n\co  \mathbf{Q}B_{n+1} \rightarrow \mathcal{M}$$ 
by assigning each element $L\in B_{n+1}$ to $S^3_{(L,1)}$.  
Note that $\tilde{\lambda}_n$ is well-defined, as $S^3_{(L,+1)}$ is an integral homology sphere 
for all $L\in B_{n+1}$.  

Denote by $\mathbf{I}$ the submodule of $\mathbf{Q}B_{n+1}$ generated by element $(L-L')$ such that 
$\tilde{\lambda}_n(L-L')$ is in $\mathcal{M}_{2n-1}$.  
The following follows immediately from \cite{hbrun} and \fullref{2n-1}.  
\begin{lem}\label{lem:lh}
Let $L$ and $L'$ be two link-homotopic (or $C^a_{n+1}$--equivalent) $(n+1)$--component Brunnian links. Then $L-L'\in \mathbf{I}$.  
\end{lem}
Note that two link-homotopic $(n+1)$--component Brunnian links satisfy $L-L'\in J_{2n+1}(n+1)$ \cite[Proposition 7.1]{hm1}.  We 
generalize \fullref{lem:lh} as follows.  
\begin{prop} \label{claim}
  Let $L,L'$ be two $(n+1)$--component Brunnian links in $S^3$ such that $L-L'\in J_{2n+1}(n+1)$.  
  Then $L-L'\in \mathbf{I}$.  
\end{prop}
\begin{proof}
Let $B$ be an $(n+1)$--component Brunnian link in $S^3$.   
By \cite[Section 5]{hm1}, we have 
$B\sim_{C^a_{n+1}} U_F$, where $F=T_1\cup ... \cup T_m$ is a simple $C^a_n$--forest $F$ for $U$ in $S^3$ such that, for all 
$1\le i\le p$, we have $T_i=T_{\sigma_i}$ for some $\sigma_i\in S_{n-1}$ 
(see \fullref{Tsigma} for the definition of $T_{\sigma_i}$).  
By \fullref{lem:lh} we thus have 
  $$ B\equiv U_{F} \textrm{ mod $\mathbf{I}$.} $$  
Observe that we have the equality 
  $$ U_{F}=\sum_{F'\subseteq F} (-1)^{|F'|} [U;F']. $$ 
For all $F'\subseteq F$, denote by $G(F')$ the clasper obtained in $S^3$ by performing $(+1)$--framed surgery along $U$.  
As in \fullref{proof2n-2}, we have $\tilde{\lambda}_n(U_F)=S^3_{(U_F,+1)}\cong S^3_{G(F)}$.  
As each $C^a_n$--tree in $F'$ is turned into a $Y_{n-1}$--tree of $S^3$ by this operation, we have 
$\tilde{\lambda}_n([U;F'])=[S^3;G(F')]\in \mathcal{M}_{(n-1).|F'|}$.  
In particular, $\tilde{\lambda}_n([U;F'])\in \mathcal{M}_{2n-2}$ for all $F'$ with $|F'|\ge 3$.  
It follows that 
  $$ B\equiv \sum_{F'\subseteq F \textrm{ / } |F'|\le 2} (-1)^{|F'|} [U;F'] \textrm{ mod $\mathbf{I}$}. $$  
By strictly the same arguments as in the proof of \cite[Theorem 7.4]{hm1}, one can check that, for every $\sigma\in S_{n-1}$, 
$[U;T_{\sigma}]\equiv \frac{1}{2}[U;T_{\sigma},\tilde{T}_{\sigma}]\textrm{ mod $\mathbf{I}$}$.  
It follows that 
  $$ B\equiv U+\frac{1}{2}\sum_{1\le i\le m}  [U;T_{\sigma_i},\tilde{T}_{\sigma_i}] 
    + \sum_{1\le i\ne j\le m}  [U;T_{\sigma_i},\tilde{T}_{\sigma_j}]\textrm{ mod $\mathbf{I}$}. $$ 
It follows that $L-L'$ is equal, modulo $\mathbf{I}$, to a linear combination of the form ($\alpha_{\sigma, \sigma'}\in \mathbf{Q}$)
 \begin{equation} \label{a}
   \sum_{\sigma, \sigma'\in S_{n-1}} \alpha_{\sigma, \sigma'} [U;T_{\sigma},\tilde{T}_{\sigma'}].    
 \end{equation}
By assumption, $L-L'\in J_{2n+1}(n+1)$.  So \eqref{a} vanishes in 
$\Br(\overline{J}_{2n}(n+1))$, and is thus mapped by $h_n^{-1}$ onto a linear combination of connected trivalent 
diagrams which vanishes in $\mathcal{A}^c_{n-1}(\emptyset)$.  
\eqref{a} is thus a linear combination of terms of the following two types.  
 \begin{enumerate}
   \item(AS) $[U;T_1,T_2]+[U;T'_1,T'_2]$, where $T_1\cup T_2$ and $T'_1\cup T'_2$ differ by the cyclic order of the three edges 
             attached to a node.  
   \item(IHX) $[U;T_1,T_2]+[U;T'_1,T'_2]+[U;T''_1,T''_2]$, where $T_1\cup T_2$, $T'_1\cup T'_2$ and $T''_1\cup T''_2$ are as 
              claspers $I$, $H$ and $X$ of \fullref{fig:ihx}.  
 \end{enumerate}
Consider a term of type $(1)$.  By \cite[Corollary 4.6]{GGP}, we have 
$\tilde{\lambda}_n([U;T_1,T_2]+[U;T'_1,T'_2])\in \mathcal{M}_{2n-1}$.  
The same holds for terms of type $(2)$ by \cite[Theorem 4.11]{GGP}.  

This completes the proof. 
\end{proof}

By \fullref{2n-2} and \fullref{claim}, we have a well-defined homomorphism 
\begin{gather*}\lambda_n\co  \Br(\overline{J}_{2n}(n+1)) \rightarrow \overline{\mathcal{M}}_{2n-2}\\
\tag*{\hbox{by setting}}
\lambda_n([L-U]_{J_{2n+1}}):=[S^3 - S^3_{(L,+1)}]_{\mathcal{M}_{2n-1}}
\end{gather*}
\subsection[Proof of \ref{homspheres}]{Proof of \fullref{homspheres}} \label{proofhomspheres}
First, we show that $\lambda_n$ actually takes its values in the connected part of the Ohtsuki filtration.  
  
Recall from \fullref{brunnianpart} that $\Br(\overline{J}_{2n}(n+1))$ is generated by elements 
$[U;T_{\sigma}\cup \tilde{T}_{\sigma'}]$, for $\sigma, \sigma'\in S_{n-1}$.  
Each component $U_i$ of $U$ intersects one disk-leaf $f_i$ of $T_{\sigma}$ and one disk-leaf $f'_i$ of $T_{\sigma'}$.  
Denote by $G_{\sigma,\sigma'}$ the $Y_{2n-2}$--graph obtained from $T_{\sigma}\cup \tilde{T}_{\sigma'}$ by connecting,  
for each $1\le i\le n+1$, the edges incident to $f_i$ and $f'_i$.  
\begin{lem} \label{inclusion}
  For all $\sigma, \sigma'\in S_{n-1}$, 
   $$ \lambda_n([U;T_{\sigma}\cup \tilde{T}_{\sigma'}])\equiv [S^3;G_{\sigma,\sigma'}] \textrm{ mod } \mathcal{M}_{2n-1}. $$ 
  Consequently, we have 
   $$ \lambda_n(\Br(\overline{J}_{2n}(n+1)))\subset \Co(\overline{\mathcal{M}}_{2n-2}). $$  
\end{lem}
\begin{proof}
 For any $\sigma, \sigma'\in S_{n-1}$, we have 
   $$ \lambda_n([U;T_{\sigma}\cup \tilde{T}_{\sigma'}]) = 
      -S^3_{G(T_{\sigma}\cup \tilde{T}_{\sigma'})}+S^3_{G(T_{\sigma})}+S^3_{G(T_{\sigma'})}-S^3, $$  
 where, if $F$ is a $C^a_n$--forest for $U$, $G(F)$ denotes the clasper obtained in $S^3$ by $(+1)$--framed surgery along $U$. 
 
 For all $\tau \in S_{n-1}$, $G(T_{\tau})$ is a linear $Y_{n-1}$--tree whose 
 leaves are all $(-1)$--special leaves.  
 So by \fullref{table}, there exists a union $G_\tau$ of $Y_k$--trees, $k\ge 2n-2$ such that 
 $ S^3_{G(T_{\tau})}\cong S^3_{G_\tau}$.  
 
 On the other hand, $G(T_{\sigma}\cup \tilde{T}_{\sigma'})$ is obtained from $T_{\sigma}\cup \tilde{T}_{\sigma'}$ 
 by replacing $f_i\cup f'_i$ by a pair of Hopf-linked $(-1)$--framed leaves (as illustrated in \fullref{surgery}), 
 for $1\le i\le n+1$.  
  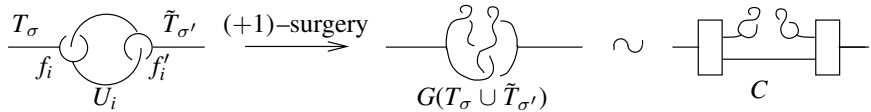
\begin{figure}[ht!]
  \bc
\begin{picture}(0,0)%
\includegraphics{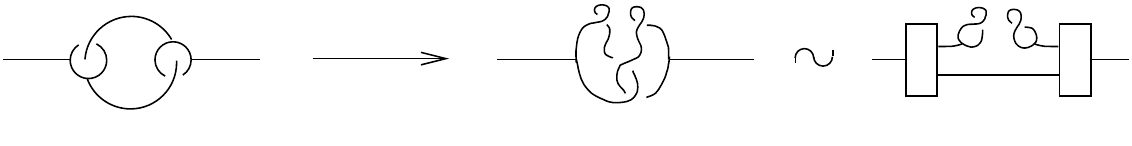}%
\end{picture}%
\setlength{\unitlength}{3947sp}%
\begingroup\makeatletter\ifx\SetFigFont\undefined%
\gdef\SetFigFont#1#2#3#4#5{%
  \reset@font\fontsize{#1}{#2pt}%
  \fontfamily{#3}\fontseries{#4}\fontshape{#5}%
  \selectfont}%
\fi\endgroup%
\begin{picture}(5432,689)(0,-1318)
\put(1351,-811){\makebox(0,0)[lb]{\smash{\SetFigFont{10}{12.0}{\rmdefault}{\mddefault}{\updefault}{\color[rgb]{0,0,0}$(+1)$--surgery}%
}}}
\put(4658,-1223){\makebox(0,0)[lb]{\smash{\SetFigFont{10}{12.0}{\rmdefault}{\mddefault}{\updefault}{\color[rgb]{0,0,0}$C$}%
}}}
\put( 23,-811){\makebox(0,0)[lb]{\smash{\SetFigFont{10}{12.0}{\rmdefault}{\mddefault}{\updefault}{\color[rgb]{0,0,0}$T_{\sigma}$}%
}}}
\put(969,-811){\makebox(0,0)[lb]{\smash{\SetFigFont{10}{12.0}{\rmdefault}{\mddefault}{\updefault}{\color[rgb]{0,0,0}$\tilde{T}_{\sigma'}$}%
}}}
\put(526,-1269){\makebox(0,0)[lb]{\smash{\SetFigFont{10}{12.0}{\rmdefault}{\mddefault}{\updefault}{\color[rgb]{0,0,0}$U_i$}%
}}}
\put(2565,-1269){\makebox(0,0)[lb]{\smash{\SetFigFont{10}{12.0}{\rmdefault}{\mddefault}{\updefault}{\color[rgb]{0,0,0}$G(T_{\sigma}\cup \tilde{T}_{\sigma'})$}%
}}}
\put(196,-1071){\makebox(0,0)[lb]{\smash{\SetFigFont{10}{12.0}{\rmdefault}{\mddefault}{\updefault}{\color[rgb]{0,0,0}$f_i$}%
}}}
\put(906,-1071){\makebox(0,0)[lb]{\smash{\SetFigFont{10}{12.0}{\rmdefault}{\mddefault}{\updefault}{\color[rgb]{0,0,0}$f'_i$}%
}}}
\end{picture}
    \caption{Performing $(+1)$--framed surgery along the unlink $U$} \label{surgery}
  \ec
  \end{figure}
 By Habiro's move 7 and 2, $G(T_{\sigma}\cup \tilde{T}_{\sigma'})$ is equivalent to the clasper $C$ obtained by replacing 
 each such pair of Hopf-linked leaves by two boxes as shown in \fullref{surgery}.  
 By using the zip construction and \fullref{crossingchange}, we obtain 
 $$ S^3_C\cong S^3_{G(T_{\sigma}\cup \tilde{T}_{\sigma'})}\sim_{Y_{2n-1}} S^3_{G_{\sigma,\sigma'}\cup G(T_{\sigma})\cup 
    G(\tilde{T}_{\sigma'})}.  $$
 It follows that
   $$ \lambda_n([U;T_{\sigma}\cup \tilde{T}_{\sigma'}])\equiv -S^3_{G_{\sigma,\sigma'}\cup G_{\sigma}\cup G_{\sigma'}} 
      + S^3_{G_{\sigma}} + S^3_{G_{\sigma'}} - S^3 \textrm{ mod } \mathcal{M}_{2n-1}. $$
  By using the equality $S^3_{G_{\sigma,\sigma'}\cup G_{\sigma}\cup G_{\sigma'}} = 
   \sum_{G'\subseteq \{ G_{\sigma,\sigma'},G_{\sigma},G_{\sigma'}\} } (-1)^{|G'|} [S^3;G']$,  
 one can easily check that 
  $$ S^3_{G_{\sigma,\sigma'}\cup G_{\sigma}\cup G_{\sigma'}}\equiv 
     S^3_{G_{\sigma,\sigma'}} + S^3_{G_{\sigma}} + S^3_{G_{\sigma'}} - 2S^3\textrm{ mod } \mathcal{M}_{2n-1}. $$ 
 \noindent (here we use the fact that $G_{\sigma,\sigma'}$ and each connected component of $G_{\sigma}$ and $G_{\sigma'}$ have 
  degree $\ge 2n-2$).  The result follows.  
\end{proof}
Clearly, the composite $\alpha_{n-1} \lambda_n$ is the map  
  $$ \kappa_n\co   \Br(\overline{J}_{2n}(n+1))\longrightarrow \overline{\mathcal{S}}_{2n-2} $$
announced in the statement of \fullref{homspheres}.  
By \fullref{e2}, it suffices to show that $\lambda_n$ is an isomorphism 
to obtain the theorem.  This is implied by the next lemma.  
\begin{lem} \label{e3}
  For $n\ge 3$, the following diagram commutes: 
    \begin{equation*}
      \xymatrix{
        \mathcal{A}^c_{n-1}(\emptyset)
        \ar[d]_{h_n}
        \ar[dr]^{\varphi_n}\\
        \Br(\overline{J}_{2n}(n+1))
        \ar[r]_{\lambda_n}&\Co(\overline{\mathcal{M}}_{2n-2})
      }
    \end{equation*}
\end{lem}
\begin{proof}
 As pointed out in \cite[Section 3.5]{hm2}, one can easily check that $\mathcal{A}^c_{n-1}(\emptyset)$ is generated by 
 the elements 
 $\Gamma_{\sigma}$ depicted in \fullref{generator}, for all $\sigma\in S_{n-1}$.  
  \begin{figure}[ht!]
  \bc
\begin{picture}(0,0)%
\includegraphics{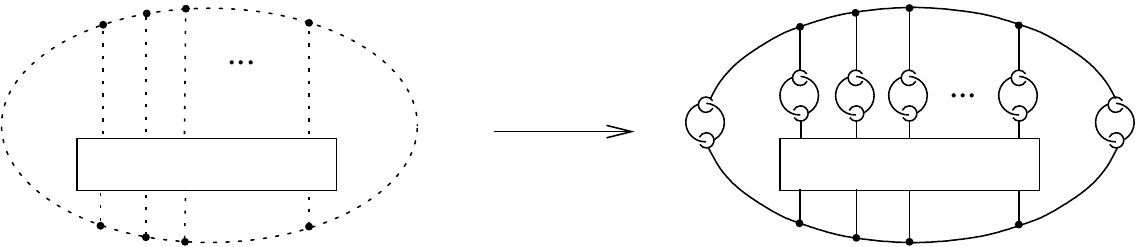}%
\end{picture}%
\setlength{\unitlength}{3947sp}%
\begingroup\makeatletter\ifx\SetFigFont\undefined%
\gdef\SetFigFont#1#2#3#4#5{%
  \reset@font\fontsize{#1}{#2pt}%
  \fontfamily{#3}\fontseries{#4}\fontshape{#5}%
  \selectfont}%
\fi\endgroup%
\begin{picture}(5453,1168)(102,-988)
\put(3938,-638){\makebox(0,0)[lb]{\smash{\SetFigFont{10}{12.0}{\rmdefault}{\mddefault}{\updefault}{\color[rgb]{0,0,0}permutation $\sigma$}%
}}}
\put(563,-638){\makebox(0,0)[lb]{\smash{\SetFigFont{10}{12.0}{\rmdefault}{\mddefault}{\updefault}{\color[rgb]{0,0,0}permutation $\sigma$}%
}}}
\put(5313,-849){\makebox(0,0)[lb]{\smash{\SetFigFont{10}{12.0}{\rmdefault}{\mddefault}{\updefault}{\color[rgb]{0,0,0}$T_{\sigma}$}%
}}}
\put(5313,-58){\makebox(0,0)[lb]{\smash{\SetFigFont{10}{12.0}{\rmdefault}{\mddefault}{\updefault}{\color[rgb]{0,0,0}$T_1$}%
}}}
\end{picture}
    \caption{The connected trivalent diagram $\Gamma_{\sigma}$, and the two simple linear 
             $C^a_{n}$--trees $T_1$ and $T_{\sigma}$} 
    \label{generator}
  \ec
  \end{figure}
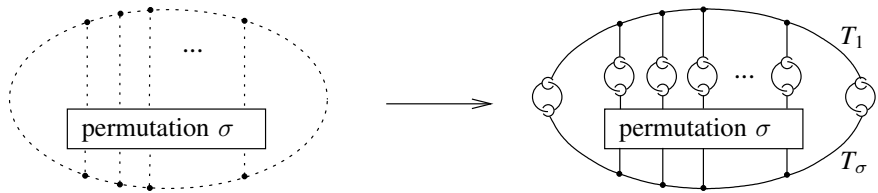
  
 For such an element $\Gamma_{\sigma}$, a representative for $h_n(\Gamma_{\sigma})$ is $[U;T_{1}\cup T_{\sigma}]$, 
 where $T_1$ and $T_{\sigma}$ are two $C^a_{n}$--trees for $U$ as represented in \fullref{generator}.  
 As seen in the proof of \fullref{inclusion}, 
 $\lambda_n([U;T_{1}\cup T_{\sigma}])=[S^3;G_{1,\sigma}]_{\mathcal{M}_{2n-1}}$, where $G_{1,\sigma}$ is 
 obtained by replacing each pair of disk-leaves intersecting the same component of $U$ by an edge.  
 Clearly, this $Y_{2n-2}$--graph satisfies $\varphi_n(\Gamma_{\sigma})=[S^3;G_{1,\sigma}]_{\mathcal{M}_{2n-1}}$.  
\end{proof}

The various results proved of this section can be summed up in the following commutative diagram ($n\ge 2$)   
  \begin{equation*}
    \xymatrix{
      &\mathcal{A}^c_{n-1}(\emptyset)
      \ar[dl]_{h_n}
      \ar[d]_{\varphi_{n-1}}
      \ar[dr]^{\phi_{n-1}}\\
      \Br(\overline{J}_{2n}(n+1))
      \ar[r]_{\lambda_n}&\Co(\overline{\mathcal{M}}_{2n-2})
      \ar[r]_{\quad \alpha_{n-1}}&\overline{\mathcal{S}}_{2n-2}
    }
  \end{equation*} 
  where all arrows are isomorphism over $\mathbf{Q}$.  
\subsection{Brunnian links with vanishing Milnor invariants}
In this last subsection, we can work over the integers.  

Habegger and Orr also studied finite type invariants of integral homology spheres obtained by $(+1)$--framed surgery along links 
in $S^3$.  In particular, \cite[Theorem 2.1]{HO} deals with $(+1)$--framed surgery along $l$--component Brunnian links with 
vanishing Milnor invariants of length $\le 2l-1$, and appears to have some similarities with our results.  

Let $\Br^{l}(\ov{J}_k(n))$ denote the $\mathbf{Z}$--submodule of $\ov{J}_k(n)$ generated by elements $[L-U]_{J_{k+1}}$ where $L$ 
is an $n$--component Brunnian link with vanishing Milnor invariants of length $\le l$.  
Let $U_{(k)}$ denote the $k$--component unlink $U_1\cup\cdots\cup U_k$ in $S^3$.  
Let 
  $$ S_{n+1}\co  \Br(\ov{J}_{2n}(n+1)) \longrightarrow \mathbf{Z}\mathcal{L}(n) $$
be the map defined by 
  $$ S_{n+1}([L-U_{(n+1)}]_{J_{2n+1}}) = s_{n+1}(L) - U_{(n)}, $$
where $s_{n+1}(L)$ denotes the $n$--component link in $S^3$ obtained by $(+1)$--framed surgery along the $(n+1)^{th}$ component of $L$.  
In particular, $s_{n+1}(U_{(n+1)})=U_{(n)}$.  

We can show that, for $n\ge 3$,   
  \begin{enumerate}
    \item $S_{n+1}(\Br(\ov{J}_{2n}(n+1)))=\Br^{2n-1}(\ov{J}_{2n-1}(n))$  
    \item $S_{n+1}\oq \co \Br(\ov{J}_{2n}(n+1))\oq \rightarrow \Br^{2n-1}(\ov{J}_{2n-1}(n))\oq $ is an isomorphism.  
  \end{enumerate}
The proof involves the same technique as in the preceding section, and makes use of Theorem 6.1 of \cite{HM}.  
\section[The proof of \ref{linear2n+1}]{The proof of \fullref{linear2n+1}} \label{zeproof}
In this section, we give the proof of \fullref{linear2n+1}.  
For that purpose, it is convenient to state a few more technical lemmas on claspers.  
\begin{lem} \label{spebox}
  The move of \fullref{box} produces equivalent claspers.  
  \begin{figure}[ht!]
  \bc
    \includegraphics{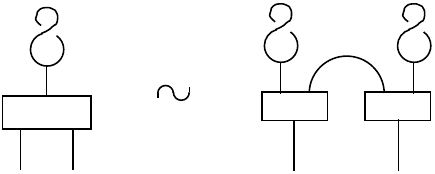}
    \caption{} \label{box}
  \ec
  \end{figure}
\end{lem} 
\noindent This is an easy consequence of \cite[Proposition 2.7]{habiro}.   
\begin{lem} \label{lemsubtree}
Let $G$ be a clasper in a $3$--manifold $M$ containing a $Y_k$--subtree $T$, $k\ge 1$, such that a branch of $T$ is incident to a box 
as shown in \fullref{F40}.  There, $e$ is an edge of $G$ which is not contained in $T$.  Then 
$$M_G\sim_{Y_{k+1}} M_{G'},$$  
where $G'$ is the clasper depicted in the right-hand side of \fullref{F40}.  
    \begin{figure}[ht!]
    \bc
    \includegraphics{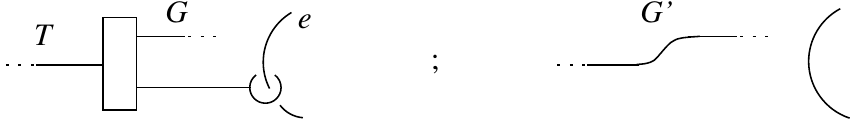} 
    \caption{} \label{F40}
    \ec
    \end{figure}
\end{lem}
\noindent The proof is omitted.  It is straightforward, and uses Habiro's move 12 and a zip construction.  
\begin{lem} \label{lemplus}
Let $G$ be a clasper in a $3$--manifold $M$ such that a $3$--ball $B$ in $M$ intersects $G$ as depicted in \fullref{F41}.  
There, the nodes $n_1$ and $n_2$ are both in a $Y_k$--subtree $T$, $k\ge 2$, and $e$ is an edge of $G$ which is not contained in $T$.  
Then 
$$M_G\sim_{Y_{k+1}} M_{G'},$$ 
where $G'$ is identical to $G$ outside of $B$, where it is as shown in \fullref{F41}.  
    \begin{figure}[ht!]
    \bc
\begin{picture}(0,0)%
\includegraphics{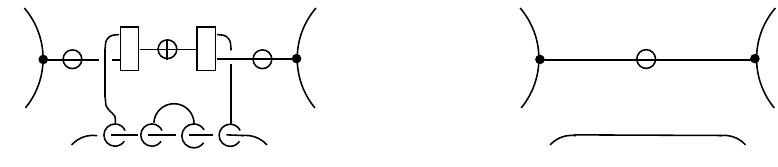}%
\end{picture}%
\setlength{\unitlength}{3947sp}%
\begingroup\makeatletter\ifx\SetFigFont\undefined%
\gdef\SetFigFont#1#2#3#4#5{%
  \reset@font\fontsize{#1}{#2pt}%
  \fontfamily{#3}\fontseries{#4}\fontshape{#5}%
  \selectfont}%
\fi\endgroup%
\begin{picture}(3730,793)(2989,-417)
\put(3701,268){\makebox(0,0)[lb]{\smash{\SetFigFont{10}{12.0}{\rmdefault}{\mddefault}{\updefault}{\color[rgb]{0,0,0}$G$}%
}}}
\put(4291,-377){\makebox(0,0)[lb]{\smash{\SetFigFont{10}{12.0}{\rmdefault}{\mddefault}{\updefault}{\color[rgb]{0,0,0}$e$}%
}}}
\put(4999,-43){\makebox(0,0)[lb]{\smash{\SetFigFont{10}{12.0}{\rmdefault}{\mddefault}{\updefault}{\color[rgb]{0,0,0};}%
}}}
\put(4456, 47){\makebox(0,0)[lb]{\smash{\SetFigFont{10}{12.0}{\rmdefault}{\mddefault}{\updefault}{\color[rgb]{0,0,0}$n_2$}%
}}}
\put(2989, 51){\makebox(0,0)[lb]{\smash{\SetFigFont{10}{12.0}{\rmdefault}{\mddefault}{\updefault}{\color[rgb]{0,0,0}$n_1$}%
}}}
\put(5985,249){\makebox(0,0)[lb]{\smash{\SetFigFont{10}{12.0}{\rmdefault}{\mddefault}{\updefault}{\color[rgb]{0,0,0}$G'$}%
}}}
\end{picture}
    \caption{} \label{F41}
    \ec
    \end{figure}
\end{lem}
\begin{proof}
 By an isotopy, $G$ is seen to be equivalent to the clasper $G_1$ represented in \fullref{F42}.  
 By applying the move of \cite[Figure 38]{habiro} to $G_1$, and then 
 applying Habiro's move 6 twice, we obtain the clasper $G_2\sim G_1$ of \fullref{F42}.  
    \begin{figure}[ht!]
    \bc
\begin{picture}(0,0)%
\includegraphics{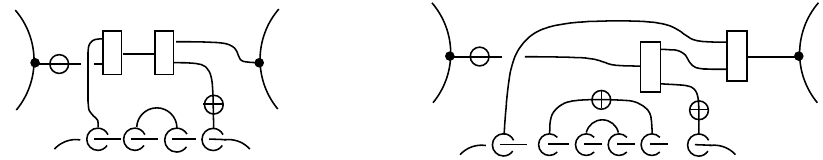}%
\end{picture}%
\setlength{\unitlength}{3947sp}%
\begingroup\makeatletter\ifx\SetFigFont\undefined%
\gdef\SetFigFont#1#2#3#4#5{%
  \reset@font\fontsize{#1}{#2pt}%
  \fontfamily{#3}\fontseries{#4}\fontshape{#5}%
  \selectfont}%
\fi\endgroup%
\begin{picture}(3937,762)(2807,-365)
\put(4757,-320){\makebox(0,0)[lb]{\smash{\SetFigFont{10}{12.0}{\rmdefault}{\mddefault}{\updefault}{\color[rgb]{0,0,0}$G_2$}%
}}}
\put(2807,-308){\makebox(0,0)[lb]{\smash{\SetFigFont{10}{12.0}{\rmdefault}{\mddefault}{\updefault}{\color[rgb]{0,0,0}$G_1$}%
}}}
\end{picture}
    \caption{} \label{F42}
    \ec
    \end{figure}
 Consider the two $\textsf{I}$--shaped claspers $I_1\cup I_2$ of $G_2$ which appear in the figure.  
 By Habiro's move 6 and 4, we have that $G_2\sim G_2\setminus (I_1\cup I_2)$.  The result then follows from \fullref{lemsubtree}.  
\end{proof}

We can now prove \fullref{linear2n+1}.  
 
Let $G$ be a linear $Y_n$--tree in a $3$--manifold $M$, $n\ge 2$, with $n+2$ $(-1)$--special leaves, and let   
$N$ denote an $s$--regular neighborhood $N$ of $G$.  As noted previously, $N$ is a $3$--ball in $M$.  

By $(n-1)$ applications of \fullref{slide}, $G$ is equivalent to the clasper $\tilde{G}$ represented in 
\fullref{F34}.  The first step of this proof is to show the following.  
\begin{claim} \label{Klaim} 
We have  
 $$ \tilde{G}\sim C, $$ 
in $N$, where $C$ is the clasper containing a $Y_{2n}$--subtree represented in \fullref{F34}.  
    \begin{figure}[ht!]
    \bc
\begin{picture}(0,0)%
\includegraphics{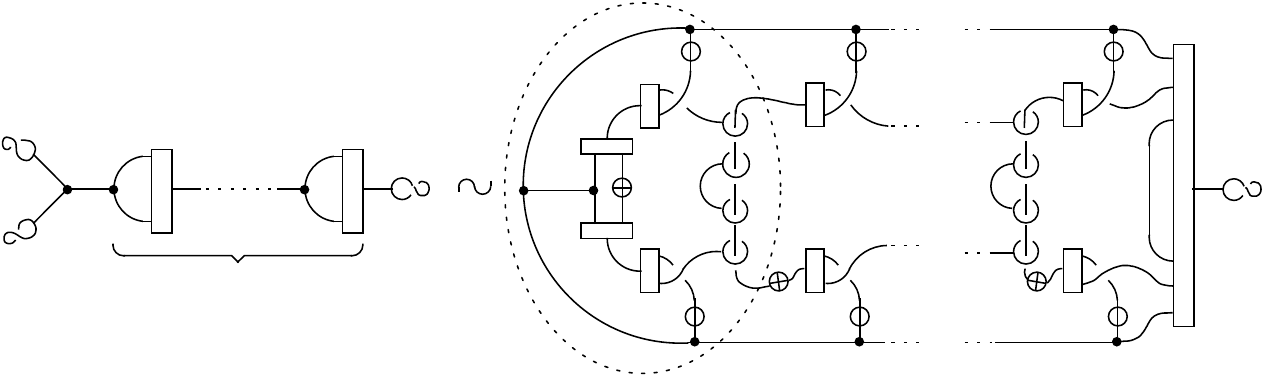}%
\end{picture}%
\setlength{\unitlength}{3947sp}%
\begingroup\makeatletter\ifx\SetFigFont\undefined%
\gdef\SetFigFont#1#2#3#4#5{%
  \reset@font\fontsize{#1}{#2pt}%
  \fontfamily{#3}\fontseries{#4}\fontshape{#5}%
  \selectfont}%
\fi\endgroup%
\begin{picture}(6064,1792)(0,-1512)
\put(826,-1111){\makebox(0,0)[lb]{\smash{\SetFigFont{10}{12.0}{\rmdefault}{\mddefault}{\updefault}{\color[rgb]{0,0,0}$(n-1)$ times}%
}}}
\put(2426,-1452){\makebox(0,0)[lb]{\smash{\SetFigFont{10}{12.0}{\rmdefault}{\mddefault}{\updefault}{\color[rgb]{0,0,0}$B$}%
}}}
\put(976,-436){\makebox(0,0)[lb]{\smash{\SetFigFont{10}{12.0}{\rmdefault}{\mddefault}{\updefault}{\color[rgb]{0,0,0}$\tilde{G}$}%
}}}
\put(4276,-661){\makebox(0,0)[lb]{\smash{\SetFigFont{10}{12.0}{\rmdefault}{\mddefault}{\updefault}{\color[rgb]{0,0,0}$C$}%
}}}
\put(1351,-586){\makebox(0,0)[lb]{\smash{\SetFigFont{10}{12.0}{\rmdefault}{\mddefault}{\updefault}{\color[rgb]{0,0,0}$v$}%
}}}
\end{picture}
    \caption{} \label{F34}
    \ec
    \end{figure}
\end{claim}
\begin{proof}
 Consider the box of $\tilde{G}$ which is connected to one $(-1)$--special leaf.  
 This box is connected to a node $v$ by two edges.  
 By applying \fullref{move} at $v$, and \fullref{lass}, we obtain the clasper represented in 
 \fullref{prooflinear} (a).  
 Then apply recursively \fullref{move} and Habiro's move 6, as shown in 
 \fullref{prooflinear} (b), until we obtain a clasper $G'\sim \tilde{G}$ with only one node 
 connected to two $(-1)$--special leaves.  See in \fullref{prooflinear} (c).  
    \begin{figure}[ht!]
    \bc
\labellist\small\hair 0pt
\pinlabel (a) [t] at 3 7 
\pinlabel (b) [t] at 70 7
\pinlabel (c) [t] at 164 7
\pinlabel (d) [t] at 267 7
\pinlabel $G'$ at 161 49
\pinlabel $G''$ <3pt, 8pt> [b] at 261 29
\pinlabel $f$ at 277 88
\pinlabel $f'$ at 328 88
\pinlabel $c$ at 269 136
\endlabellist
    \includegraphics{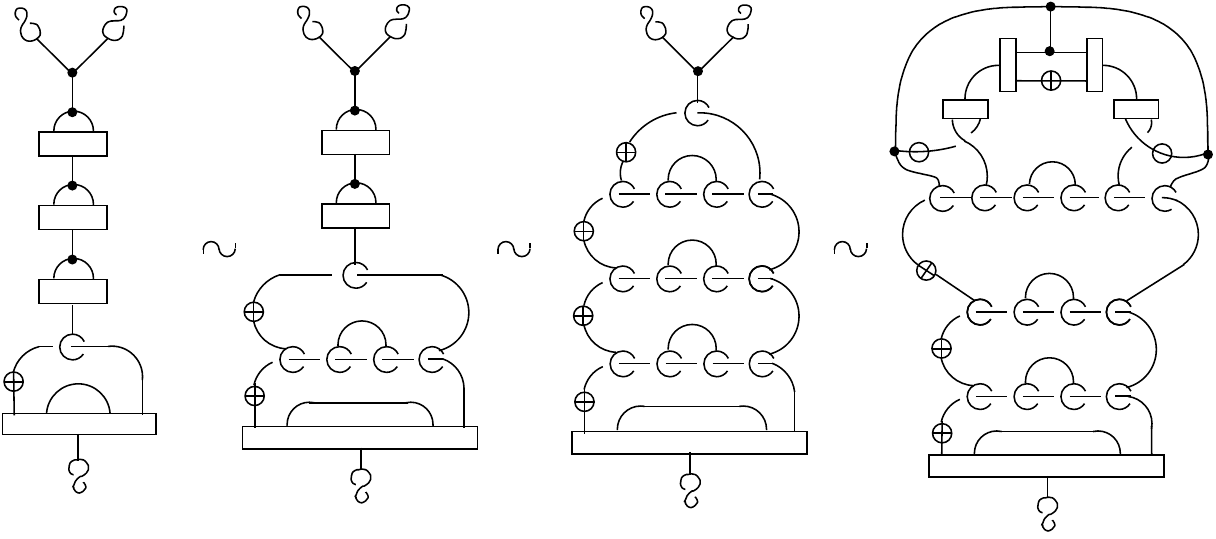} 
    \caption{Here, for simplicity, we consider the case $n=5$.} \label{prooflinear}
    \ec
    \end{figure}
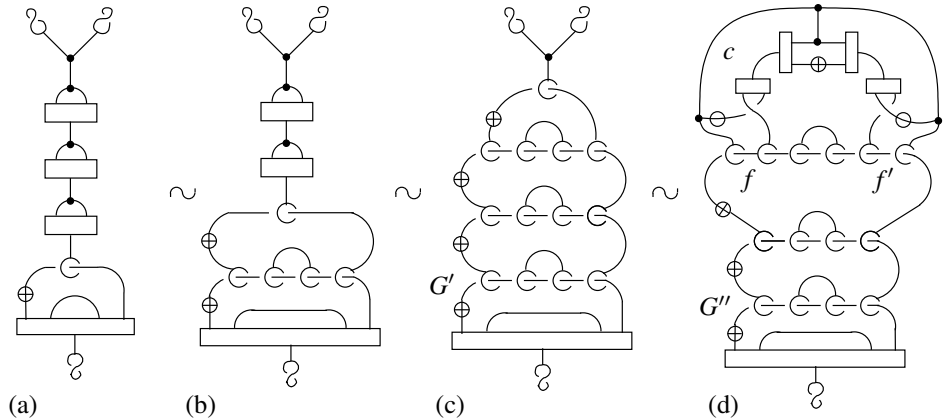 
 By applying the move of \fullref{prooftwo} and Habiro's move 6, we have 
 $G'\sim G''$, where $G''$ contains a component $c$ with $4$ nodes and with two leaves $f$ and $f'$ lacing an edge $e$ -- see \fullref{prooflinear} (d).\footnote{   
 Here we say that a leaf of a clasper $G$ \emph{laces} an edge if it forms an unknot which bounds a disk $D$ with respect to which 
 it is $0$--framed, such that the interior of $D$ intersects $G$ once, transversally, at an edge.  }  
 We can apply Habiro's move 12 to these two leaves, and then Habiro's move 6 to create two 
 new leaves lacing an edge.  Apply recursively these two moves until no new leaf lacing an edge is 
 created: the result is the desired clasper $C$ which contains a $Y_{2n}$--subtree $T$, as represented in \fullref{F34}.  
\end{proof}

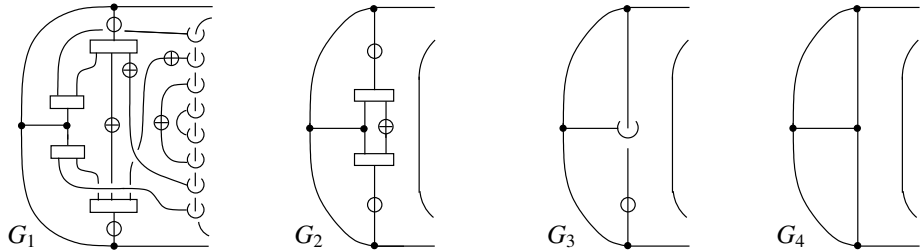
\begin{figure}[ht!]
    \bc
\begin{picture}(0,0)%
\includegraphics{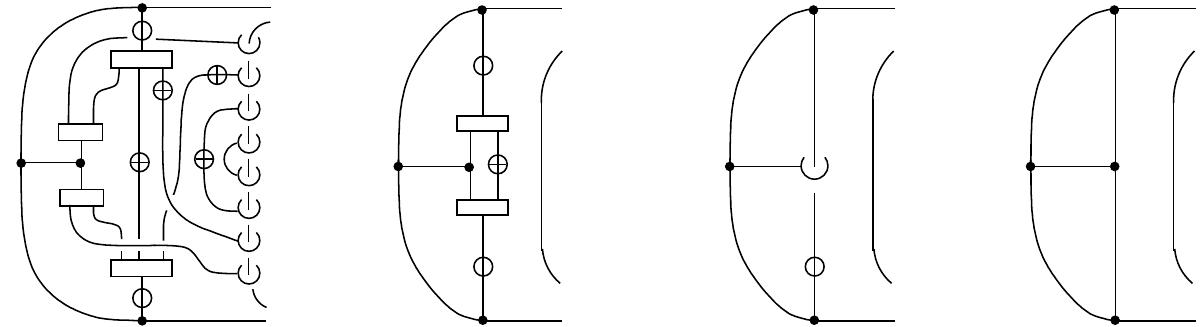}%
\end{picture}%
\setlength{\unitlength}{3947sp}%
\begingroup\makeatletter\ifx\SetFigFont\undefined%
\gdef\SetFigFont#1#2#3#4#5{%
  \reset@font\fontsize{#1}{#2pt}%
  \fontfamily{#3}\fontseries{#4}\fontshape{#5}%
  \selectfont}%
\fi\endgroup%
\begin{picture}(5754,1556)(29,-1383)
\put( 29,-1337){\makebox(0,0)[lb]{\smash{\SetFigFont{10}{12.0}{\rmdefault}{\mddefault}{\updefault}{\color[rgb]{0,0,0}$G_1$}%
}}}
\put(3418,-1337){\makebox(0,0)[lb]{\smash{\SetFigFont{10}{12.0}{\rmdefault}{\mddefault}{\updefault}{\color[rgb]{0,0,0}$G_3$}%
}}}
\put(1831,-1337){\makebox(0,0)[lb]{\smash{\SetFigFont{10}{12.0}{\rmdefault}{\mddefault}{\updefault}{\color[rgb]{0,0,0}$G_2$}%
}}}
\put(4859,-1337){\makebox(0,0)[lb]{\smash{\SetFigFont{10}{12.0}{\rmdefault}{\mddefault}{\updefault}{\color[rgb]{0,0,0}$G_4$}%
}}}
\end{picture}
    \caption{These four claspers are identical to $C$ outside $B$.} \label{F36}
    \ec
    \end{figure}

Consider in $N$ a $3$--ball $B$ which intersects $C$ as depicted.  
By several applications of the move of \cite[Figure 38]{habiro} and of Habiro's move 6, we obtain the clasper 
$G_1\sim G$ which is identical to $C$ outside $B$, where it is as shown in \fullref{F36}.  
By Habiro's move 6 and 4, we can freely remove the pair of $\textsf{I}$--shaped claspers which appear in the figure 
(see the proof of \fullref{lemplus}).  By further applying four times \fullref{lemsubtree}, we thus obtain 
the clasper $G_2$ of \fullref{F36}, which satisfies $N_{G_2}\sim_{Y_{2n+1}} N_{G_1}$.  
By an isotopy, we can apply Habiro's move 12 to show that $N_{G_2}\sim N_{G_3}$, 
where $G_3$ is as shown in \fullref{F36}.  
By using \cite[page 398]{ohtsuki}, we obtain $N_{G_3}\sim_{Y_{2n+1}} N_{G_4}$.\footnote{
We use the up-most figure of \cite[page 398]{ohtsuki}.  The arguments given there are 
for graph claspers, but they can be used in our situation.  }

\begin{figure}[ht!]
    \bc
\begin{picture}(0,0)%
\includegraphics{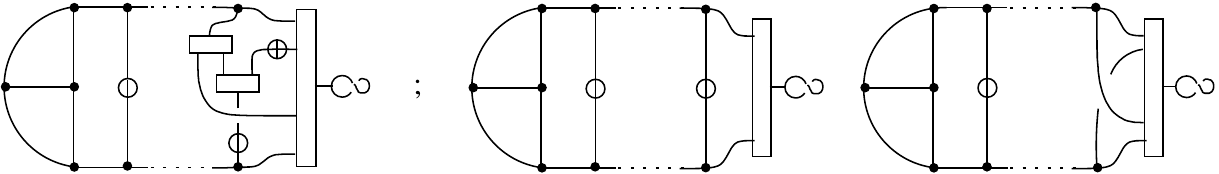}%
\end{picture}%
\setlength{\unitlength}{3947sp}%
\begingroup\makeatletter\ifx\SetFigFont\undefined%
\gdef\SetFigFont#1#2#3#4#5{%
  \reset@font\fontsize{#1}{#2pt}%
  \fontfamily{#3}\fontseries{#4}\fontshape{#5}%
  \selectfont}%
\fi\endgroup%
\begin{picture}(5838,824)(14,-1857)
\put(4866,-1746){\makebox(0,0)[lb]{\smash{\SetFigFont{10}{12.0}{\rmdefault}{\mddefault}{\updefault}{\color[rgb]{0,0,0}$G_7$}%
}}}
\put(756,-1746){\makebox(0,0)[lb]{\smash{\SetFigFont{10}{12.0}{\rmdefault}{\mddefault}{\updefault}{\color[rgb]{0,0,0}$G_5$}%
}}}
\put(3001,-1746){\makebox(0,0)[lb]{\smash{\SetFigFont{10}{12.0}{\rmdefault}{\mddefault}{\updefault}{\color[rgb]{0,0,0}$G_6$}%
}}}
\end{picture}
    \caption{} \label{F35}
    \ec
    \end{figure}

Observe that $G_4$ satisfies the hypothesis of \fullref{lemplus}.  Actually, we can apply 
\fullref{lemplus} recursively ($n-3$) times.  By further applying, to the resulting clasper, 
strictly the same arguments as in the proof of \fullref{lemplus}, we obtain 
$N_{G_4}\sim_{Y_{2n+1}} N_{G_5}$, where $G_5$ is the clasper shown in \fullref{F35}.  
It follows, by the zip construction and \fullref{crossingchange}, that 
 $$ N_{G_5}\sim_{Y_{2n+1}} N_{G_6\cup G_7}, $$ 
where $G_6$ and $G_7$ are two disjoint claspers in $N$ as represented in \fullref{F35}.  
    
By \fullref{spebox} and \fullref{table} (for $l=1$), it is not hard to check that $N_{G_7}\sim_{Y_{2n+1}} N$  
and that $N_{G_6}\sim_{Y_{2n+1}} N_{\Theta_n}$.  

This concludes the proof of \fullref{linear2n+1}.  

\bibliographystyle{gtart}
\bibliography{link}

\end{document}